\documentclass[11pt]{amsart}

\usepackage{amssymb,amsthm,amsmath,latexsym,mathrsfs}
\usepackage{amsfonts}
\usepackage[numbers,sort&compress]{natbib}
\usepackage{color}
\usepackage{graphicx}
\usepackage[colorlinks,
            linkcolor=blue,
            anchorcolor=blue,
            citecolor=blue
            ]{hyperref}

\usepackage{amsfonts}
\usepackage{latexsym,amssymb,mathrsfs}

\let\newpf\proof \let\proof\relax 
\newenvironment{pf}{\newpf[\proofname]}{\qed\endtrivlist}

\newcommand{\ba}{\overline{A}}

\def\be{\begin{equation}}
\def\ee{\end{equation}}

\def\ba{{\begin{align}}}
\def\ea{{\end{align}}}

\def\bm{\begin{matrix}}
\def\em{\end{matrix}}

\def\0{{\mathbf 0}}

\newtheorem{Theorem}{Theorem}[section]
\newtheorem{Lemma}{Lemma}[section]
\newtheorem{Proposition}{Proposition}[section]

\newtheorem{Definition}{Definition}[section]

\numberwithin{equation}{section}

\def \bn {\hfill \\ \smallskip\noindent}

\theoremstyle{definition}

\def\proof{\bn {\bf Proof.} }

\def\ssm{\smallsetminus}

\renewcommand{\setminus}{\ssm}

\newcommand{\C}{{\mathbb C}}

\newcommand{\R}{{\mathbb R}}
\newcommand{\T}{{\mathbb T}}

\newcommand{\Z}{{\mathbb Z}}

\def\B0{{\bold{0}}}

\catcode`\@=12

\def\Empty{}
\newcommand\oplabel[1]{
  \def\OpArg{#1} \ifx \OpArg\Empty {} \else
    \label{#1}
  \fi}

\newcommand{\comm}[1]{}
\newcommand{\comment}[1]{}

\begin{document}

\title[CS]{Explicit construction of Gevrey Quasi-Periodic Discrete Schr\"odinger Operators with
cantor spectrum}

\author{Xuanji  Hou,
 Li Zhang}

\date{\today}

\setcounter{tocdepth}{1}

\begin{abstract}
  We  construct  1-dim difference Schr\"odinger operators with a class of Gevrey potentials such that Cantor spectrum occurs together with the estimations of open spectral gaps for $0<|\lambda|\leq 1$.  The proof is based on  KAM  and Moser-P\"oschel argument .
\end{abstract}

\maketitle

\section{INTRODUCTION AND MAIN RESULTS}
The difference Schr\"odinger operators  on $l^{2}(\mathbb{Z})$. of  the form
\begin{equation}\label{1.1}
  (H_{\lambda v,\alpha,\theta}u)_{n}=u_{n+1}+u_{n-1}+ \lambda v(\theta+n\alpha)u_{n},\quad \forall ~n\in\mathbb{Z},
\end{equation}
are concerned with  in this article,  where $\theta\in\mathbb{T}^{d}:=(\mathbb{R}/2\pi\mathbb{Z})^{d}$ is called the {\it phase},~$v:\mathbb{T}^{d}\rightarrow\mathbb{R}$ is  called the {\it potential}, $0<|\lambda|\leq1 $~and $\alpha\in\mathbb{T}^{d}$ is  called the {\it frequencies}.~When $(1,\alpha)$ is rationally independent\footnote{We say that $(1,\alpha)$ is \textit{rationally independent} if  $\langle k,\alpha\rangle+j\neq 0$ for all $(k,j)\in\Z^{d+1}$.}, the spectrum of \eqref{1.1} is a bounded close set of $\mathbb{R}$ independent on the phase,  and we then use $\Sigma_{\lambda v,\alpha}$ to denote it.   We say that \eqref{1.1} has {\it Cantor spectrum} if $\Sigma_{\lambda v,\alpha}$ is a Cantor set, i.e.,  $\R\setminus \Sigma_{\lambda v,\alpha}$ is a dense open set.

Observe that $\R\setminus \Sigma_{\lambda v,\alpha}$ is open and then is the union of countable open intervals which are called {\it spectral gaps}.  Denote by $N_{\lambda v,\alpha}(E)$ the  {\it integrated density of states} (IDS) of $H_{\lambda v,\alpha,\theta}$, given as
 $N_{\lambda v,\alpha}(E):=\int_{\mathbb{T}}\mu_{\lambda v,\alpha,\theta}(-\infty,E]d\theta$,
 where $\mu_{ \lambda v,\alpha,\theta}$ is the spectral measure.  The \textit{Gap-Labelling Theorem} \cite{JM82} states that all spectral gaps can be formulated as
 $$I_k:= N_{\lambda v,\alpha}^{-1}(\langle k, \alpha\rangle)=\{E\in\R\,|\, N_{\lambda v,\alpha}(E)=\langle k, \alpha\rangle\}$$
 for some $k\in \Z^d$.  However, conversely,  for any $k\in\Z^d$, $I_k$ is a spectral gap if an only if $|I_k|>0$, and we prefer to to call it a {\it candidate  spectral gap}.


\subsection{Cantor Spectrum Of Quasi-Periodic Schr\"odinger Operators }
As is well - known, there are many conclusions regarding Cantor spectrum  for  quasi-periodic Schr\"odinger operators.  Avila-Bochi-Damanik \cite{ABD09} proved that Cantor spectrum holds generally in $C^{0}$-sense for any fixed rational independent frequencies. Eliasson \cite{E92} proved  that for given Diophantine frequencies $\alpha\in {\rm DC}_d(\gamma,\tau)$,~1-dim  differential Schr\"odinger operators  have Cantor spectrum for  generic small analytic potential by Moser-P\"oschel argument \cite{MP84} (the smallness of the potential depend of the analyticity and $\gamma,\tau$).  Here we say that $\alpha\in\mathbb{R}^{d}$ is {\it Diophantine } if there exist $\gamma>0, \tau>d-1$ such that $\alpha\in {\rm DC}_{d}(\gamma,\tau)$,~where
\begin{equation}\label{1.2}
  {\rm DC}_{d}(\gamma,\tau):=\{x\in\mathbb{R}^{d}:\inf_{j\in\mathbb{Z}}|\langle n,x\rangle-2\pi j|>\frac{\gamma}{|n|^{\tau}},\quad \forall~n\in\mathbb{Z}^{d}\backslash\{0\}\}.
\end{equation}
The Eliasson's result can be extended to difference Schr\"odinger operators and also $C^{\infty}$ , $C^k$ regularities \cite{Am09,CCYZ17,CG17}.  Define ${\rm DC}_{d}:=\cup_{\gamma>0,\tau>d-1}{\rm DC}_{d}(\gamma,\tau)$. When $d=1$ and $\alpha\in  {\rm DC}_{1}$,   non-perturbative version of Eliasson's  result  also holds, ~i.e.,~the smallness of the potential is independent of the frequency \cite{P06}.  However,  one cannot obtain any explicit example of Cantor spectrum from these results.

The Ten Martini problem conjectures that the Cantor spectrum holds for  \eqref{1.1} with $d=1$ and $v(\theta)=2cos2\pi\theta$ (called {\it almost Mathieu operator}), which is well-known long-standing problem and was solved thoroughly by Avila-Jitomirskaya's renowned work \cite{AJ08} using several techniques.  One can refer to Avila-You-Zhou \cite{AYZ} for the proof of {\it Dry Ten Martini Problem} furthermore asks whether for any $\lambda\neq0$ and irrational $\alpha$,~all candidate  spectral gaps  are non-collapsed, and  the noncritical case ($\lambda\neq 1$) .

There are also apart from Cantor spectrum examples besides almost Mathieu. Sinai and Wang-Zhang \cite{S87,WZ}  established Cantor spectrum for large $C^{2}$  cosine-type potentials. Recently, Cantor spectrum examples for operators close to almost Mathieu are  given  by Ge-Jitomirskaya-You \cite{GJY23}. ~In \cite{HSY}, as for  1-dim differential Schr\"odinger operators, ~Hou-You-Shan presented  a strategy  to explicitly construct a class of  small {\it Gevrey} quasi-periodic potential such that Cantor spectrum holds using only the information of $\alpha$. Similar construction for  1-dim difference Schr\"odinger operators with $C^k$ potential was given by He-Cheng  \cite{HC23}. In this article,~we  construct discrete  1-dim dierence Schr\"odinger operators with a class of Geverey potentials such that Cantor spectrum occurs.




\begin{Theorem}\label{main-th}
Given $s\in(0,\frac{1}{2})$ and  $\alpha\in {\rm DC}_d(\gamma,\tau)$.~One can construct explicitly a set $\mathcal{K}\subseteq\mathbb{Z}^{d}$ depending on $\alpha$ and $s$,~such that for the Schr\"odinger operator $H_{\lambda v,\alpha,\theta}$ \eqref{1.1} with
\begin{equation}\label{1.3}
  v(\theta)=\Sigma_{k\in\mathcal{K}}e^{-|k|^{s}}cos\langle k,\theta\rangle,
\end{equation}
is a Gevrey real function, $\Sigma_{\lambda v,\alpha}$ is Cantor for all $0<|\lambda|\leq1$.
\end{Theorem}

The function $Hu=Eu$, for $H=H_{\lambda v,\alpha,\theta}$ \eqref{1.1} can be transformed to
\begin{equation}
\label{E-cocycle}
(\alpha, S_{E}^{\lambda v}):\quad \T^d\times\R^2\rightarrow \T^d\times\R^2,\quad (\theta, x)\mapsto (\theta+\alpha, S_{E}^{\lambda v}(\theta)x),
\end{equation}
which is called a {\it cocycle}, where $S_{E}^{\lambda v}=\left(
             \begin{array}{cc}
               E-\lambda v(\theta) & -1\\
               1 & 0 \\
             \end{array}
           \right)$.  Define
$$S_{E}^{(\lambda v,n)}:=S_{E}^{\lambda v}(\cdot+(n-1)\alpha)\cdots S_{E}^{\lambda v}(\cdot+\alpha)S_{E}^{\lambda v},\qquad n=1,2,\cdots.$$
A is well-known, $E\in \R\backslash \Sigma_{\lambda v,\alpha}$ as long as $(\alpha,S_{E}^{\lambda v})$ is {\it uniformly hyperbolic}, i.e.,
 $\|S_{E}^{(\lambda v,n)}\|\geq c \zeta^n$ $(\zeta>1,c>0)$ holds for all $n$.
To prove the Cantor spectrum, we need to demonstrate that  $\{E\in\R\,|\,\hbox{$(\alpha,S_{E}^{\lambda v})$ is uniformly hyperbolic.}\}$ is dense. The cocycle  $(\alpha,S_{E}^{\lambda v})$ is said to be {\it reducible} if  it can be conjugated to some constant cocycle, i.e,
\begin{equation}
\label{const.-E-cocycle}
(\alpha, C):\quad \T^d\times\R^2\rightarrow \T^d\times\R^2,\quad (\theta, x)\mapsto (\theta+\alpha,Cx)
\end{equation}
where $C\in SL(2,\R)$ is some constant matrix. Our aim is to prove that, for $v$  given in Theorem \ref{main-th} and for $E$ in a dense set of $\R$,  $(\alpha,S_{E}^{\lambda v})$ can be conjugated to \eqref{const.-E-cocycle} with $C$ being a hyperbolic matrix (the eigenvalues  are not on $\{z\in\C\,|\,|z|=1\}$), which implies that $(\alpha,S_{E}^{\lambda v})$ is uniformly hyperbolic.

One fundamental tool  to the problem of reducibility is  KAM  \cite{DS75, MP84, E92}  etc..~KAM is an iteration method and the main difficulties come from the  {\it small divisors problems} . The KAM scheme  in this article is somewhat different from the usual one. More precisely, at each  KAM step, besides the usual KAM estimations, we need further estimations related to the special form of $v$, and in the end we obtain the reducibility with the reduced constant matrix being hyperbolic or parabolic (eigenvalues are multiple $1$), provided that
\begin{equation}\label{1.4}
  \rho(E)\in\{\frac{1}{2}\langle k,\alpha\rangle|k\in\mathcal{K}\}.
\end{equation}
Then, Moser-P\"oschel can be applied to ensure that all  gaps with $k\in \mathcal{K}$ are open, and the structure of $\mathcal{K}$ ensure the Cantor spectrum (density of spectral gaps).  The crucial point is to design $\mathcal{K}$, which is similar to the one in \cite{HSY}.

\subsection{Estimates On Spectral Gap}
The  estimates on the spectral gaps are also of  great interests to us.  The upper bound estimates play an important role in proving spectrum homogeneity \cite{SP95,SP97},~which is a crucial subject in the study of inverse spectral theory.  There are some results on the  estimation  of the spectral gaps given by  Amor \cite{HA09}, Damanik-Goldstein \cite{DGo14} and etc.. In Leguil-You-Zhao-Zhou \cite{LYZZ}, there introducs the upper bounds estimations of spectral gaps for Shr\"odinger operator with small quasi-periodic analytic potentials with Diophantine frequencies, and gives both upper bounds and lower bounds estimations of spectral gaps for the almost Mathieu operator. In this article, as for the  quasi-periodic Shr\"odinger operator with Diophantine frequencies and the potential $ \lambda v$ defined by  \eqref{1.3} with $0<|\lambda|\leq 1$,  one can get both upper bounds and lower bounds estimations for chosen spectral gaps which are dense on $\R$.

\begin{Theorem}\label{th-1.2}
Let $\alpha\in {\rm DC}_d(\gamma,\tau)$.~
For all $k\in \mathcal{K}$,   the candidate spectral gap $I_{k}(\lambda v)=(E_{k}^{-},E_{k}^{+})$ of $H_{\lambda v,\alpha,\theta}$ satisfies
\begin{equation}\label{1.9}
  |\lambda|^{2} e^{-\frac{13}{6}|k|^{2s}}\leq E_{k}^{+}- E_{k}^{-} \leq \sqrt{|\lambda|}e^{-\frac{3}{20}|k|^{s}}.
\end{equation}
\end{Theorem}

\par
\subsection{Structure of the paper}
Let us provide an introduction to the remaining parts of this article.~In Section 2,~we introduce some necessary notations and lemmata and prove that the existence of the sets $\mathcal{K}$ .~In Section 3,~we prove one step of KAM,~with some delicate estimates,~based on some abstract KAM lemmata and we utilize Proposition \ref{lem3.2} to construct a KAM iteration,~with crucial estimates.~In Section 4,~we prove Theorem \ref{th-1.2} by applying Lemma \ref{lem4.1}.~In Section 5,~Theorem \ref{th-6.1} provides both the upper and lower bounds for any label $k\in\mathcal{K}$.

\section{PRELIMINARIES }

Let  $M_2(\R)$($M_2(\C)$)  denote the space of all $2\times 2$ real (complex) matrices.\footnote{$M_2(\R)$ ($M_2(\C)$) is equipped
the usual operator norm $\|\cdot\|$ induced by Euclidean norm of $\R^2$ (unitary norm of $\C^2$).} Denote by $SL(2,\R)$  the  group $\big\{\left(
\begin{array}{ccc}
   a &  b\cr
   c &  d \end{array}
   \right)\,\big|\,ad-bc=1,\, a,b,c,d\in \R\big\}$  with  $sl(2,\R)=\big\{\left(
\begin{array}{ccc}
   u &  v\cr
   w &  - u \end{array}
   \right)\,\big|\,u,v,w\in \R\big\}$ being its Lie algebra, and
by  $SU(1,1)$  the  group  $\big\{\left(
\begin{array}{ccc}
   a &  b\cr
   \bar{b} &  \bar{a} \end{array}
   \right)\,\big|\,|a|^2-|b|^2=1,\, a,b\in \C\big\}$ with $su(1,1)=\big\{\left(
\begin{array}{ccc}
   i\rho &  w\cr
   \overline{w} &  - i\rho \end{array}
   \right)\,\big|\,\rho\in\R, w\in \C\big\}$ being its  Lie algebra.

 $SL(2,\R)$ ,~ $(sl(2,\R))$ ,~ $SU(1,1)$ and  $(su(1,1))$ are isomorphic via algebraic
 conjugation through  $P_{\lozenge}= \frac{1}{\sqrt{-2i}}\left(
\begin{array}{ccc}
   -i & -1\cr
   -i & 1 \end{array}
   \right)$.  Let   $PSL(2,\R)$ and $PSU(1,1)$ denotes quotient groups $SL(2,\R)/\{\pm I\}$ and
   $SU(1,1)/\{\pm I\}$ respectively.

\begin{Definition}\label{Def-trc}
For any $A\in SL(2,\R)$ (or $A\in SU(1,1)$), we define
\[|A|_{trc}\triangleq |c|\] provided that
 $UAU^*=\left(
\begin{array}{ccc}
   * &  c\cr
   0 &   * \end{array}
   \right)$ for some $U\in U(2)$ and $c\in\C$.
\end{Definition}

The $|\cdot|_{trc}$  is  well-defined.~For any  $A\in SL(2,\R)$ (or $A\in SU(1,1)$) with
 $spec(A)=\{\pm\mu\}$ ($\mu\in\C$),
 \begin{equation}\label{est.-trc}
|A|_{trc}\leq \|A\|\leq |A|_{trc}+|\mu|.
\end{equation}

\begin{Lemma}\cite[Lemma 3.2]{HSY}  \label{Lem-trc-basic}
For any  $A_1,A_2\in SL(2,\R)$ (or $A_1,A_2\in SU(1,1)$),
\begin{equation}
|A_2|_{trc}\leq |A_1|_{trc}+ 2\|A_1-A_2\|.
\end{equation}
\end{Lemma}

\begin{Lemma}\cite[Lemma 3.3]{HSY}\label{Lem-algebra-diag}
Let $A\in SL(2,\R)$ ($A\in SU(1,1)$) with $spec(A)=\{e^{\pm i\rho}\}\,(\rho\neq 0)$. There is $P\in SU(1,1)$,
such that $PAP^{-1}=\left(
\begin{array}{ccc}
   e^{i\rho} &  0\cr
   0 &   e^{-i\rho}\end{array}
   \right)$, and
\begin{eqnarray}
\|P\|^2\leq 2|\rho|^{-1}\|A\| \leq  2(1+|\rho|^{-1} )|A|_{trc}.
\end{eqnarray}
\end{Lemma}

\begin{Lemma}\label{Lem-A1}
For $C\in su(1,1)$ with $spec(C)=\{\pm\mu\}$, and
\begin{equation}\label{A1-trans}
  C=\left(
      \begin{array}{cc}
        ia & b \\
        \bar{b} & -ia \\
      \end{array}
    \right), \quad( a\in\mathbb{R}, b\in\mathbb{C}),
\end{equation}
then $A\triangleq e^{C}\in SU(1,1)$ with
\begin{equation}\label{A2-trans}
A=\left(
    \begin{array}{cc}
      cosh(\mu)+ia\frac{sinh(\mu)}{\mu} & b \frac{sinh(\mu)}{\mu}\\
      \bar{b}\frac{sinh(\mu)}{\mu} & cosh(\mu)-ia\frac{sinh(\mu)}{\mu} \\
    \end{array}
  \right).
\end{equation}
\end{Lemma}
\begin{pf}
By simple calculation, it follows.
\end{pf}

\begin{Lemma}\label{Lem-A2}
Let $A=M^{-1}\left(
         \begin{array}{cc}
           a & b \\
           \bar{b} & \bar{a} \\
         \end{array}
       \right)M \in SL(2,\mathbb{R})$ with spec$(A)={e^{\pm i\rho}}$, where $|a|^{2}-|b|^{2}=1$ and $a,b\in\mathbb{C}$, then exists $R_{\phi}:=\left(\begin{array}{cc}
cos2\pi\phi & -sin2\pi\phi\\
sin2\pi\phi & cos2\pi\phi \\
\end{array}
\right)$
for any $\phi\in\mathbb{R}$, satisfying $R_{-\phi}AR_{\phi}=\left(
                                                              \begin{array}{cc}
                                                                1 & |b| \\
                                                                0 & 1 \\
                                                              \end{array}
                                                            \right)$.

\end{Lemma}
\begin{pf}
Note that $|a|^{2}-|b|^{2}=1$, we have $\Im(a)=|b|$. Thus $MAM^{-1}=\left(
                                                                      \begin{array}{cc}
                                                                        1+i|b| & b \\
                                                                        \bar{b} & 1-i|b| \\
                                                                      \end{array}
                                                                    \right)$, the desired result can be reached with some computations.
\end{pf}

\subsection{Quasi-periodic~Cocycles~And~Fibered~Rotation~Number}
  Given $A\in C({\T}^{d},SL(2,\mathbb{C}))$ and rationally independent $\alpha\in\mathbb{R}^{d}$,~we define the quasi-periodic $cocyle(\alpha,A)$:
 \begin{equation*}
 	\begin{split}
 	(\alpha,A):\quad \mathbb{T}^{d}\times\mathbb{C}^{2} \rightarrow\mathbb{T}^{d}\times\mathbb{C}^{2},\quad
 		(\theta,v)  \mapsto(\theta+\alpha,A(\theta)v).\\
 	\end{split}
 \end{equation*}
The iterates of $(\alpha,A)$ are of the form $$(\alpha,A)^{\circ n}=(\alpha,A)\circ\cdots \circ(\alpha,A)=(n\alpha,\mathcal{A}_{n}),$$
  where
\begin{equation*}
  	\mathcal{A}_{n}:=
  	\begin{cases}
  		A(\cdot+(n-1)\alpha)\cdots A(\cdot),& n>0,\\
  		I,&n=0,\\
  		A(\cdot-n\alpha)^{-1}\cdots A(\cdot-\alpha) ^{-1}&n<0.
  	\end{cases}
\end{equation*}
$(\alpha,A)$ is called  $uniformly ~hyperbolic$ if,~for every $x\in\mathbb{T}^{d}$,~there exists a continuous splitting $\mathbb{C}^{2}=E^{s}\bigoplus E^{u}$ such that for every $n\geq 0$,
$$\Vert\mathcal{A}_{n}(\theta)^{-1}\omega\Vert\leq C\lambda^{n}\Vert\omega\Vert,\quad \omega\in E^{s}(\theta),$$
$$\Vert\mathcal{A}_{n}(\theta-n\alpha)^{-1}v\Vert\leq C\lambda^{n}\Vert v\Vert,\quad v\in E^{u}(\theta),$$
for some constants $C,~c>0$.~This splitting is invariant in the sense that
$$A(\theta)E_{s}(\theta)=E_{s}(\theta+\alpha),\qquad A(\theta)E_{u}(\theta)=E_{u}(\theta+\alpha).$$

When  $A\in C({\T}^{d},SL(2,\mathbb{C}))$ is homotopic to the identity, ~$(\alpha, A)$ induces the projective skew-product
$F_{A}:\mathbb{T}^{d}\times\mathbb{S}^{1}\rightarrow\mathbb{T}^{d}\times\mathbb{S}^{1}$ with
$$F_{A}(\theta,w):=(\theta+\alpha,\frac{A(\theta)\cdot w}{|A(\theta)\cdot w|}),$$
which is also homotopic to the identity and can be lifted to a map $\tilde{F}_{A}:\mathbb{T}^{d}\times\mathbb{R}\rightarrow\mathbb{T}^{d}\times\mathbb{R}$ of the form $\tilde{F}_{A}(\theta,y)=(\theta+\alpha,y+\psi_{\theta}(y))$,  and the $fibered~rotation~number$ of $(\alpha,A)$ is
given as  $$\rho(\alpha,A):=\int_{\mathbb{T}^{d}\times\mathbb{R}}\psi_{\theta}(y))d\mu(\theta,y)\quad( mod\quad\mathbb{Z}),$$~where $\mu$ is an invariant probability measure on
$\mathbb{T}^{d}\times\mathbb{R}$, and the fibered rotations number does not depend on the lift and the measure $\mu$.  One can refer to \cite{JM82} for more details.

Define  $R_{\phi}:=\left(\begin{array}{cc}
cos2\pi\phi & -sin2\pi\phi\\
sin2\pi\phi & cos2\pi\phi \\
\end{array}
\right)$
for any $\phi\in\mathbb{R}$. We obviously have the following conclusion:

\begin{Lemma}\label{lem2.1}
If $A\in C(\mathbb{T}^{d},SL(2,\mathbb{R}))$ is homotopic to the identity,~then
$$|\rho(\alpha,A)-\phi|<|A-R_{\phi}|_{\mathbb{T}^{d}}.$$
for all $\phi\in\R$.
\end{Lemma}

If $A:\mathbb{T}^{d}\rightarrow PSL(2,\mathbb{R})$ is homotopic to $\theta\mapsto R_{\frac{\langle k,\theta\rangle}{2}}$ for some $k\in\mathbb{Z}^{d}$,~then we call $k$ is the $degree$ of $A$ and denote it by deg $A$.~

\begin{Lemma}\label{lem2.2}
If the cocycle $(\alpha,A_{1})$ is conjugated to the cocycle  $(\alpha,A_{2})$,~i.e.,~$B(\theta+\alpha)^{-1}A_{1}(\theta)B(\theta)=A_{2}(\theta)$,~for some
$B:\mathbb{T}^{d}\rightarrow PSL(2,\mathbb{R})$. When $B\equiv P\in SL(2,\mathbb{R})$ is constant or $B$ is $C^{1}$ smooth and homotopic to $\theta\mapsto I$, we have
$$\rho(\alpha,A_{1})=\rho(\alpha,A_{2}).$$
When $B(\theta)=\left(
                  \begin{array}{cc}
                    cos \frac{\langle k,\theta\rangle}{2} & -sin \frac{\langle k,\theta\rangle}{2} \\
                    sin \frac{\langle k,\theta\rangle}{2} & cos \frac{\langle k,\theta\rangle}{2} \\
                  \end{array}
                \right)
$ for some $k\in\mathbb{Z}^{d}$
$$\rho(\alpha,A_{1})=\rho(\alpha,A_{2})+\frac{\langle k,\alpha\rangle}{2}\qquad (mod\quad\Z).$$
\end{Lemma}

A typical example is given by  $Schr\ddot{o}dinger$ cocycles $(\alpha,S_{E}^{\lambda v})$,~with
$$S_{E}^{\lambda v}(\cdot):=\left(
                      \begin{array}{cc}
                        E-\lambda v(\cdot) & -1 \\
                        1 & 0 \\
                      \end{array}
                    \right),\quad E\in\mathbb{R},
$$
which  were introduced due to the connection with the eigenvalue equation $H_{\lambda v,\alpha,\theta}\mu=E\mu$.~There are close relationships between the spectral set of $H_{\lambda v,\alpha,\theta}$ and the dynamics of $(\alpha,S_{E}^{\lambda v})$,  as indicated by  the well-known fact:~
$E\in\Sigma_{\lambda v,\alpha}$ if and only if $(\alpha,S_{E}^{\lambda v})$ is not uniformly hyperbolic.


\subsection{On Functions}
Let $*$ denotes $\R$, $\C$ or a set of matrix.  For any integrable $*$ valued function
$F$ on the $d$-dimensional torus $\T^d=\R^d/2\pi\Z^d$,
\begin{equation}
F(\theta)\thicksim \sum_{k\in\Z^d}\widehat{F}(k)e^{i\langle k,\theta \rangle},\quad \widehat{F}(k)\triangleq\frac{1}{(2\pi)^d}\oint_{\T^d}F(\theta)e^{-i\langle k,\theta \rangle}d\theta,
\end{equation}
and we also use the notation
\begin{equation}
\langle F \rangle \triangleq\widehat{F}(0).
\end{equation}
In particular, we denote by
\begin{equation}
\widehat{F}(k)_{p,q} \qquad (k\in\Z^d, \,1\leq p, q\leq m )
\end{equation}
the $(p,q)$ entry of $\widehat{F}(k)$ for  $m\times m$ matrix  valued  integrable function $F$.

Given $h\geq 0$, we introduce the \textit{Wiener norm} of $F$ as
\begin{equation}
|F|_h\triangleq \sum_{k\in\Z^d}\|\widehat{F}(k)\|e^{|k|h},
\end{equation}
where  $\|\cdot\|$ is absolute value, complex number modular, matrix norm , as the $*$ may be. We use $\mathcal{B}_{h}(\T^d, * )$  to denote the set of all
  $F: \T^d \rightarrow *$ with $|F|_h<\infty$. Obviously,
\begin{eqnarray}
\mathcal{B}_{h}(\T^d, * )\subseteq \mathcal{B}_{h_+}(\T^d, * )\quad (h>h_+\geq 0).
\end{eqnarray}

  For any $N>0$,
 the truncation operators $\mathcal{T}_N$ and $\mathcal{R}_N$  are given as
 \begin{eqnarray}
\qquad\qquad &&\mathcal{T}_N F=\sum_{k\in\Z^d,\,|k|< N}\widehat{F}(k)e^{i\langle k,\theta \rangle}, \quad
 \mathcal{R}_N F=\sum_{k\in\Z^d,\,|k|\geq N}\widehat{F}(k)e^{i\langle k,\theta \rangle}
 \end{eqnarray}
for $F\in \mathcal{B}_{h}(\T^d, * )$. For any $F\in \mathcal{B}_{h}(\T^d, * )$ $(h\geq 0)$,
\begin{eqnarray}
\qquad
&& F= \mathcal{T}_N F +\mathcal{R}_N F,\qquad |F|_h=|\mathcal{T}_N F|_h +|\mathcal{R}_N F|_h; \\
&&|F|_{h_+}\leq |F|_h,\quad
|\mathcal{R}_N F|_{h_+}\leq |F|_he^{-N(h-h+)}, \quad  \forall~ h_+\in[0,h].
\end{eqnarray}

Let us  also define that, for any $N>0$,
\begin{equation}
\mathcal{B}_{h}^{(\geq N)}(\T^d, * )\triangleq \mathcal{R}_N \mathcal{B}_{h}(\T^d, * )=\{
\mathcal{R}_N F\,|\, F \in  \mathcal{B}_{h}(\T^d, * ) \}.
\end{equation}

\subsection{On Algebraic Conjugations}
Consider estimations of  $BWB^{-1}$, where $W\in \mathcal{B}_h(\T^d, su(1,1))$, $B\in \mathcal{B}_h(\T^d, PSU(1,1))$.

\begin{Lemma}\cite[Lemma 3.4]{HSY}\label{lem2.5}
Let  $A\in SU(1,1)$ with $spec(A)=\{e^{\pm i \rho}\}\neq \{0\}$, $P\in SU(1,1)$ satisfying $PAP^{-1}=\left(
\begin{array}{ccc}
   e^{i\rho }&  0\cr
   0 &   e^{-i\rho} \end{array}
   \right)$,  $W=\left(
\begin{array}{ccc}
   i u  &   w\cr
   \overline{w}& -iu\end{array}
   \right)\in \mathcal{B}_h(\T^d, su(1,1))$, and write $PWP^{-1}=\left(
\begin{array}{ccc}
   i u_+  &   w_+\cr
   \overline{w}_+ & -iu_+\end{array}
   \right)$. Then,
\begin{eqnarray}\label{W-P-alg-k-0}
\quad |\langle w_+ \rangle|&\geq &\frac{1}{2} |\rho|^{-1}\|P\|^{-2}\|[A, \langle W \rangle]\|,
\end{eqnarray}
and for all $k\in\Z^d { \backslash} \{0\}$,
\begin{equation}\label{W-P-alg-k-neq-0}
\qquad |\hat{w}_+ (k) |\geq   \frac{\|P\|^2+1}{2}(|\hat{w} (k)|-3\max\{|\hat{w}(-k)|,\,|\hat{u}(k)|\}).
\end{equation}
\end{Lemma}

\begin{Lemma}\cite[Lemma 3.5]{HSY} \label{lem2.6}
Let  $B\in \mathcal{B}_h(\T^d, PSU(1,1))$  satisfying $|B-I|_{h}\leq \frac{1}{2}$, $W\in \mathcal{B}_h(\T^d, su(1,1))$. Then,
\begin{equation}\label{est-W-conj-id}
|BWB^{-1} -W |_h \leq 4 |B-I|_{h}|W|_{h},
\end{equation}
with $W=\left(
\begin{array}{ccc}
   i u  &   w\cr
   \overline{w}& -iu\end{array}
   \right)$ and $BWB^{-1}=\left(
\begin{array}{ccc}
   i u_+  &   w_+\cr
   \overline{w}_+ & -iu_+\end{array}
   \right)$,
\begin{equation}\label{est-w-k-conj-id}
|\hat{w}_+ (k) | \geq |\hat{w} (k) |-4 |B-I|_{h}|W|_{h}e^{-|k|h},\quad k\in\Z^d.
\end{equation}
\end{Lemma}

\subsection{Continued Fraction Expansion}
Assume that $\alpha\in\mathbb{R}\backslash\mathbb{Q}$ with the denominators of best rational approximations $(q_{n})_{n\in\mathbb{N}}$ and the sequence $(N_{j})_{j\in\mathbb{N}}$ is the one defined by \eqref{k-2.8}. Given any $\alpha\in(0,1)\setminus \mathbb{Q}$, we define
$$a_{0}=0,\quad \alpha_{0}=\alpha,$$
and inductively for $k\geq1$,
$$a_{k}=[\alpha_{k-1}^{-1}],\quad \alpha_{k}=\alpha_{k-1}^{-1}-a_{k},$$
where $[\alpha]:=max\{m\in\mathbb{Z}:m\leq\alpha\}$.

Define  $p_{0}=0,~p_{1}=1,~q_{0}=1,~q_{1}=a_{1}$,~and inductively,
$$p_{k}=a_{k}p_{k-1}+p_{k-2},~\quad q_{k}=a_{k}q_{k-1}+q_{k-2}.$$
There are estimations:
$$\|k\alpha\|_{\mathbb{R}/\mathbb{Z}}\geq\|q_{n-1}\alpha\|_{\mathbb{T}},\quad for ~1\leq k<q_{n},$$
and
\begin{equation}\label{k-2.6}
  \frac{1}{q_{n}-q_{n+1}}<\|q_{n}\alpha\|_{\mathbb{R}/\mathbb{Z}}\leq\frac{1}{q_{n+1}},
\end{equation}
where $\|x\|_{\mathbb{R}/\mathbb{Z}}:=\inf_{p\in\mathbb{Z}}|x-p|$.

\subsection{Construction Of $\mathcal{K}$}
Set $s\in(0,\frac{1}{2})$ and let $N_{*}\in\mathbb{N}$ with
\begin{equation}\label{k-2.7}
  N\geq N_{*}\triangleq max\{200^{\frac{1}{1-2s}},~e^{\frac{100}{s^{2}}},~ln(|\alpha|+1),~\gamma,~\tau\},
\end{equation}
and denote
\begin{equation}\label{k-2.8}
  N_{j}=N^{\frac{12}{s}+j-1},\quad j=1,2,\cdots
\end{equation}
The set $\mathcal{K}\in\mathbb{Z}^{d}$~is chosen such that
\begin{equation}\label{k-2.9}
  \overline{\{\frac{1}{2}\langle k,\alpha\rangle|k\in\mathcal{K}\}}\footnote{$\overline{\{\cdot\}}$ is the closure set of $\{\cdot\}$.}=\mathbb{R},
\end{equation}

\begin{equation}\label{k-2.10}
  \sharp\{k\in\mathcal{K}|N_{j}\leq|k|<N_{j+2}\}\leq1,\quad j=1,2,\cdots.
\end{equation}
So, we have
\begin{equation}\label{k-2.11}
  \sharp\{k\in\mathcal{K}|\frac{21N_{j}}{10}\leq|k|<N_{j+1}\}=0,
\end{equation}
and
\begin{equation}\label{k-2.12}
  \{k\in\mathcal{K}||k|<N_{*}\}=\emptyset.
\end{equation}

\begin{Lemma}\label{k-lem-2.2}
As $\alpha$ is described above, then there exists $q_{n_{j_{*}}}\in(q_{n})_{n\in\mathbb{N}}$ with
\begin{equation}\label{k-2.13}
q_{n_{j_{*}}}\in[\frac{21N_{j}}{20},~\frac{41N_{j}}{20}]
\end{equation}
such that
\begin{equation}\label{k-2.14}
  \|q_{n_{j_{*}}}\alpha\|_{\mathbb{R}/\mathbb{Z}}:=\min_{p\in\mathbb{Z}}|q_{n_{j_{*}}}\alpha-p|<3q_{n_{j}}^{-1},
\end{equation}
where,~for the fixed $j\in\mathbb{N}$,~$q_{n_{j}}$ is the one such that $q_{n_{j}}<N_{j}\leq q_{n_{j+1}}$.
\end{Lemma}

 The proof is  essentially encompassed in Lemma 2.7 of \cite{HC23},~so we can  construct $\mathcal{K}\in\mathbb{Z}^{d}$ satisfying \eqref{k-2.9}-\eqref{k-2.12}.~Then,~ we can estimate the approximation of rational numbers to irrational number $\alpha\in\mathbb{R}\backslash\mathbb{Q}$ with  an auxiliary lemma 2.8 of \cite{HC23}.
\section{KAM PROPOSITION}
It's well-know that Theorem \ref{main-th} follows from Theorem \ref{th-1.2}.~In this section,~we construct a sequence of changes of variables which conjugate the cocycle $(\alpha,S_{E}^{v})$ to  a sequence of systems converging to a constant system.
\subsection{One Step Of KAM}
To investigate the reducibility of cocycle,~we look for near congruence $B=e^{Y(\theta)}$ to make the original cocycle $(\alpha,Ae^{F})$ to the normal system $(\alpha,C)$,~i.e
$$Ad(e^{Y(\theta)})(\alpha,Ae^{F})=(\alpha,C).$$
The linear homology equation corresponding to the above formula is
\begin{equation}\label{3.01}
  Y-A^{-1}Y(\cdot+\alpha)A=F.
\end{equation}
To find out $Y(\theta)$ in formula \eqref{3.01},~we introduce an operator:
$$\mathcal{A}: C_{h}^{\omega}(\mathbb{T}^{d},su(1,1))\rightarrow C_{h}^{\omega}(\mathbb{T}^{d},su(1,1)).$$
Obviously,~if the operator $\mathcal{A}^{-1}$ is bounded,~ it is known from the implicit function theorem that formula \eqref{3.01} has a solution.~But by  investigating  the Fourier expansion of $Y(\theta)$,~the small denominator that appears leads to \eqref{3.01} not to be solved directly.~In order to solve this obstacle,~we  adopt the following treatment.

To decompose the space $C_{h}^{\omega}(\mathbb{T}^{d},su(1,1))$,~$C_{h}^{\omega}(\mathbb{T}^{d},su(1,1))$ means that for any $h>0,~A\in SU(1,1),~\eta>0$,~we decompose the $Banach algebra$~$\mathcal{B}_{h}=C_{h}^{\omega}(\mathbb{T}^{d},su(1,1))=\mathcal{B}_{h}^{nre}(\eta)\oplus
\mathcal{B}_{h}^{re}(\eta)$, where non-resonant subspace $\mathcal{B}_{h}^{nre}(\eta)$ contains Fourier components satisfying the non-resonance conditions: $F^{nre}(\theta)=\sum_{k\in \Lambda_{1}\cup \Lambda_{2}} \widehat{F}(k)e^{i\langle k,\theta \rangle}$ and resonant subspace $\mathcal{B}_{h}^{re}(\eta)$ contains remaining Fourier components: $F^{re}(\theta)=\sum_{k \notin \Lambda_{1}\cup \Lambda_{2}} \widehat{F}(k)e^{i\langle k,\theta \rangle}$ with
\begin{equation}\label{fj1}
  \Lambda_{1}=k\in\mathbb{Z}^{d}:|\langle k, \omega \rangle|\geq\eta,\quad \Lambda_{2}=k\in\mathbb{Z}^{d}:|2\rho\pm\langle k, \omega \rangle|\geq\eta.
\end{equation}
Moreover,~we set $\mathbb{P}_{nre}(F)=F^{(nre)}$ and $\mathbb{P}_{re}(F)=F^{(re)}$ ,~respectively,~using truncation operators and the exponential decay of Fourier coefficients, we derive: $|\mathbb{P}_{re}(F)|_{h}\leq C\eta^{-1}|F|_{h}$ ensuring the stability of the decomposition.

In the KAM iteration, the non-resonant terms $\mathbb{P}_{nre}(F)$ are eliminated via homological equations, while the resonant terms $\mathbb{P}_{re}(F)$ are addressed by parameter adjustments or Floquet theory. For instance, in the elliptic case, a rotation matrix $Q(\theta)$ transforms the system into a simpler form, which is further reduced using Floquet theory.
 So the Banach space $\mathcal{B}_{h}=C_{h}^{\omega}(\mathbb{T}^{d},su(1,1))$ is decomposed into a direct sum of non-resonant and resonant subspaces, such that for any $Y\in \mathcal{B}_{h}^{nre}(\eta)$ there is
 $$A^{-1}Y(\theta+\alpha)A\in \mathcal{B}_{h}^{nre}(\eta), \quad |A^{-1}Y(\theta+\alpha)A-Y(\theta)|_{h}\geq\eta|Y(\theta)|_{r},$$
This decomposition underpins the KAM iteration, enabling the reducibility analysis of perturbed systems.

Based on the above decomposition,~we arrive at the following conclusion:
\begin{Proposition}\label{lem3.1}
Assume that $A\in SU(1,1)$ and $\eta\leq(2\|A\|)^{-4}$.~Then,~for any $F\in C_{h}^{\omega}(\mathbb{T},su(1,1))$ with $|F|_{h}<\eta^{\frac{9}{4}}$,~there exists $Y\in\mathcal{B}^{nre}(\eta)$,~$\tilde{F}\in \mathcal{B}^{re}(\eta)$ such that
    \begin{equation}\label{c-3.1}
    e^{Y(\theta+\alpha)}(Ae^{F(\theta)})e^{-Y(\theta)}=Ae^{\tilde{F}^{(re)}(\theta)},
    \end{equation}
    with estimates
     \begin{equation}\label{c-3.2}
     |Y|_{h}\leq2\eta^{-1}|F|_{h},\quad |\tilde{F}^{(re)}-\mathbb{P}_{re}F|_{h}\leq2\eta^{-7}|F|^{2}_{h}.
     \end{equation}
\end{Proposition}

\begin{pf}
The proof  is essential contained in Lemma 3.5 of \cite{HC23};~however,~through calculation we find that result
\begin{equation*}
  |\tilde{F}^{(re)}-\mathbb{P}_{re}F|_{h}\leq2\eta^{-7}|F|^{2}_{h}
\end{equation*}
is also true.
\end{pf}

\begin{Proposition}\label{lem3.2}
Let $\alpha\in\ DC(\kappa,\tau),\kappa,r>0,\tau>d-1.$ ~If $A\in SU(2,\mathbb{R}),~F\in C^{\omega}_{h}(\mathbb{T}^{d},su(2,\mathbb{R})).$~
Then for any $h^{'}\in(0,h),$~there is $c=c(\kappa,\tau,d)$ and a constant $D$ such that if
    \begin{equation}\label{3.2}
    |F|_{h}\leq\epsilon\leq\frac{c}{\Vert A\Vert^{D}}(h-h^{'})^{D\tau},
    \end{equation}
    then there exist $B\in C^{\omega}_{h^{'}}(\mathbb{T}^{d},~PSU(2,\mathbb{R})),~A_{+}\in SU(2,\mathbb{R}),~
    F_{+}\in C^{\omega}_{h^{'}}(\mathbb{T}^{d},su(2,\mathbb{R}))$ such that
    $$B^{-1}(\theta+\alpha)(Ae^{F(\theta)})B(\theta)=A_{+}e^{F_{+}(\theta)}.$$
	Let $N=\frac{2}{h-h^{'}}|\ln\epsilon|.$~Then we have the following:\\
	\textbf{(Non-resonant case:)} Assume that for $n\in\mathbb{Z}^{d},~0<|n|\leq N,$~we have $$|\langle n,\alpha\rangle|\geq\epsilon^{\frac{1}{10}},\quad
    |2\rho-\langle n,\alpha\rangle|\geq\epsilon^{\frac{1}{10}},$$
	then
    \begin{equation}\label{3.3}
    |B-Id|_{h^{'}}\leq\epsilon^{\frac{1}{2}},\quad |f_{+}|_{h^{'}}\leq4\epsilon^{2},
    \end{equation}
    \begin{equation}\label{3.4}
     |A|_{trc}\leq\frac{1}{4}\epsilon^{-\frac{1}{10}},\quad \|Ae^{\langle F\rangle}-A\|\leq2\Vert A\Vert{\epsilon}.
    \end{equation}
	\textbf{(Resonant case:)} For $n\in\mathbb{Z}^{d},~~0<|n|\leq \widetilde{N},~~\widetilde{N}\geq2N,$~if there exists $n_{*},~0<|n_{*}|\leq N$ such that $$|\langle n,\alpha\rangle|\geq\epsilon^{\frac{1}{10}},\quad
    |2\rho-\langle n,\alpha\rangle|\geq\epsilon^{\frac{1}{10}},\quad n\neq n_{*}$$
    $$|2\rho-\langle n_{*},\alpha\rangle|<\epsilon^{\frac{1}{10}},$$
	then
    \begin{equation}\label{3.5}
    |B|_{r^{'}}\leq\epsilon^{-\frac{1}{1600}}\cdot\epsilon^{-\frac{r^{'}}{r-r^{'}}},\quad |F_{+}|\ll\epsilon^{1600}.
    \end{equation}
	With the estimate
    \begin{equation}\label{3.6}
    A_{+}=\widetilde{A}^{'}e^{H}=e^{A^{''}},~~\|A^{''}\|\leq2\epsilon^{\frac{1}{10}},\quad
    |b_{+}|\leq20\epsilon^{-\frac{1}{10}}|F|_{h}e^{-|n_{*}|h'}\leq\frac{1}{100}\epsilon^{\frac{1}{10}}.
    \end{equation}
Where
$A_{+}=\left(
           \begin{array}{cc}
             ia_{+} & b_{+} \\
             \bar{b}_{+} & -ia_{+} \\
           \end{array}
         \right)
$ and $F=\left(
       \begin{array}{cc}
         if & g \\
         \bar{g} & -if \\
       \end{array}
     \right)$.
\end{Proposition}

Given that this theorem has been extensively discussed in the existing literature \cite{HY12}, to maintain the coherence of the text, we have organized the complete proof process in the appendix. Readers interested in it can refer to the appendix for detailed derivation details.

\subsection{KAM Iteration Lemma}

The quasi-periodic cocycle defined in \eqref{E-cocycle} can be rewritten as
\begin{equation}\label{n-3.7}
  \left(
    \begin{array}{c}
      u_{n+1} \\
      u_{n} \\
    \end{array}
  \right)
=(A_{E}+F_{0}(\theta+n\alpha))\left(
                                              \begin{array}{c}
                                                u_{n} \\
                                                u_{n-1} \\
                                              \end{array}
                                            \right)
\end{equation}
with
$$A_{E}=\left(
          \begin{array}{cc}
            E & -1 \\
            1 & 0 \\
          \end{array}
        \right),\quad F_{0}=\left(
                              \begin{array}{cc}
                                -\lambda v & 0 \\
                                0 & 0 \\
                              \end{array}
                            \right).
$$
Furthermore,~take $\{k_{j}\}_{j\in\mathbb{N}}\in\mathcal{K}\subset\mathbb{Z}^{d}$
$$v(\theta)=\Sigma_{j\in\mathbb{N}}e^{-|k_{j}|^{s}}cos(\langle k_{j},\theta\rangle)$$
If we denote
\begin{equation}\label{n-3.8}
  \widetilde{W}=\left(
                  \begin{array}{cc}
                    0 & 0 \\
                    1 & 0 \\
                  \end{array}
                \right),\quad W_{E}=M\widetilde{W}M^{-1},
\end{equation}
then we also have
\begin{equation}\label{n-3.9}
  S_{E}^{\lambda v}(\theta):=\left(
                       \begin{array}{cc}
                         E-\lambda v(\theta) & -1\\
                         1 &0 \\
                       \end{array}
                     \right)=A_{E}+F_{0}=A_{E}e^{\widetilde{F}}=A_{E}\prod_{k_{j}\in\mathcal{K}}e^{\lambda v_{j}(\theta)\widetilde{W}},
\end{equation}
where
$$v_{j}(\theta)=e^{-|k_{j}|^{s}}cos(\langle k_{j},\theta\rangle),\quad \widetilde{F}=\left(
                                                                                       \begin{array}{cc}
                                                                                         0 & 0 \\
                                                                                         \lambda v & 0 \\
                                                                                       \end{array}
                                                                                     \right).
$$
and $A_{E}, W_{E}\in su(1,1)$ satisfying
\begin{equation}\label{est-AEWE}
  |A_{E}|_{trc}\leq 1, \quad \|W_{E}\|\leq 1,\quad \|A_{E}\|\|W_{E}\|\leq 1,\quad \|[A_{E},W_{E}]\|=1.
\end{equation}

Let us introduce some notations(recalling $N$ given in \eqref{k-2.7}):
\begin{eqnarray}
  \label{2}
  &&h_{j}=\frac{1}{10}(N_{j+1})^{s-1},~~~j=0,1,\cdots; \\
  \label{3}
 && \mathcal{Z}_{j}\triangleq\{k\in\mathbb{Z}^{d}|N_{j}\leq|k|<N_{j+1}\},~~~j=1,2,\cdots.
\end{eqnarray}
Note that by \eqref{k-2.10} we have
\begin{equation}\label{4}
  \sharp(\mathcal{K}\cap(\mathcal{Z}_{j}\cup\mathcal{Z}_{j+1}))\leq1,~~~j=1,2,\cdots.
\end{equation}
with $v=\Sigma_{k\in\mathcal{K}}e^{-|k|^{s}}cos\langle k,\theta\rangle$ is then of the $Cocyle~(\alpha,Ae^{\Sigma_{j=1}^{\infty}\lambda v_{j}(\theta)W_{E}})$
where
\begin{equation}\label{5}
  v_{j}(\theta)=\left\{
                  \begin{array}{ll}
                    e^{-|k_{j}|^{s}}cos\{\langle k_{j},\theta\rangle\}, & if ~\{k_{j}\}\triangleq\mathcal{K}\cap\mathcal{Z}_{j}; \\
                   0, & if~\mathcal{K}\cap\mathcal{Z}_{j}=\emptyset
                  \end{array}
                \right.
\end{equation}
It is easy to see that
\begin{equation}\label{6}
  |\hat{v}(\pm k_{j})|=\frac{1}{2}e^{-|k_{j}|^{s}},~~~|v_{j}|_{h_{j}}\leq e^{-\frac{9}{10}|k_{j}|^{s}}.
\end{equation}
By the assumption that $\alpha$ is $Diophantine$ and the definition of $N$ in \eqref{k-2.7},~for all $j\geq1$ we have
\begin{equation}\label{8}
  \gamma(500N^{2s}N_{j+2})^{\gamma}(|\alpha|+1)\leq e^{\frac{1}{3000}N^{s-1}N_{j}^{s}},
\end{equation}
which ensure that,~for $0<|k|\leq 40N^{s}N_{j+2}$,
\begin{equation}\label{9}
  |\langle k,\alpha\rangle|\geq e^{-\frac{1}{50}N_{j}^{s}},
\end{equation}
and for any $\varrho\in\mathbb{R}$,~there is at most one $k$ with $0<|k|\leq 200N^{s}N_{j+1}$,~such that
\begin{equation}\label{10}
  |2\varrho-\langle k,\alpha\rangle|<e^{-\frac{1}{50}N_{j}^{s}}.
\end{equation}
Moreover,~by the definition of $N$ in \eqref{k-2.7},~for all $j\geq1$ we have
\begin{equation}\label{11}
  \begin{array}{ll}
    e^{\frac{1}{2000}N^{s-1}N_{j}^{s}}\geq max\{10^{5},~40N^{s}N_{j+1}(|\alpha|+1)\}, \\
    N_{j}^{s}\geq max \{200^{\frac{s}{1-2s}},~10^{\frac{5s}{1-s}},~200j\},~~~\frac{N_{j+1}}{N_{j}}=N\geq200^{\frac{1}{s}}.
  \end{array}
\end{equation}

  Denote by $\mathcal{NR}(L,\delta)$ all the matrices in $su(1,1)$ such that whose eigenvalues $\varrho$ satisfy
$|2\varrho-\langle k,\omega\rangle|\geq\delta$ for all any $k\in\{k\in\mathbb{Z}^{d}|0<|k|<L\}$,~and the complementary set of $\mathcal{NR}(L,\delta)$ is denoted by $\mathcal{RS}(L,\delta)$.
Starting from the system $Cocyle ~(\alpha,A_{j}e^{F_{j}+\lambda Ad(\widetilde{B}_{j}).(\Sigma_{p=j}^{\infty}v_{p}(\theta)W_{E})})$ and repeatedly applying Proposition \ref{lem3.1} or Proposition \ref{lem3.2} will lead to the following conclusion.

\begin{Lemma}\label{lem4.1}
	$\forall$ $0<h_{j}<h_{0},\gamma>0,\tau>d-1$, and $\alpha\in DC_{d}(\gamma,\tau)$. Consider the cocycle $(\alpha, A_{j}e^{\tilde{F}_{j}(\theta)})$  where $A_{j}\in SL(2,\mathbb{R})$ and
\begin{equation}\label{def-F}
\tilde{F}_{j}\triangleq F_{j}+\lambda Ad(\widetilde{B}_{j}).(v_{j}W_{E})
\end{equation}
then exist $B_{j}\in C^{\omega}_{h_{j}}(\mathbb{T}^{d},PSL(2,\mathbb{R})),$ such that $$B_{j}(\theta+\alpha)A_{j}e^{\tilde{F}_{j}(\theta)}B_{j}^{-1}(\theta)=A_{j+1}e^{F_{j+1}(\theta)},$$
with follow estimates,
\begin{equation}\label{4.1}
|F_{j}|_{\frac{3}{4}{h}_{j-1}}\leq |\lambda|\|W_{E}\| e^{-N_{j+1}^{s}},\quad |\widetilde{B}_{j}|_{{h}_{j}}\leq e^{\frac{1}{40}N_{j}^{s}},
\end{equation}
and the estimate (\ref{4.1}) with $j+1$ in place of $j$. Moreover, the following conclusions also hold:\\
(a) When $A_{j}\in\mathcal{NR}(40N^{s}N_{j+1},e^{-\frac{1}{50}N_{j}^{s}}) $ we have
\begin{equation}\label{4.15}
|B_{j}-I|_{h_{j}}\leq|\lambda|\|W_{E}\| e^{-\frac{3}{5}N_{j}^{s}},\|A_{j}-A_{j+1}\|\leq |\lambda|\|W_{E}\| e^{-\frac{3}{5}N_{j}^{s}}
\end{equation}
\begin{equation}\label{4.14}
\left\{
    \begin{array}{ll}
      \|A_{j+1}-I\|\geq\|A_{j}e^{\langle \lambda Ad(\widetilde{B}_{j}).(v_{j}W_{E})\rangle}-I\|-|\lambda|\|W_{E}\| e^{-\frac{7}{5}|k_{j}|^{s}} , & \mathcal{K}\cap \mathcal{Z}_{j}\neq\emptyset; \\
      \|A_{j+1}-I\|\geq\|A_{j}-I\|-|\lambda|\|W_{E}\| e^{-\frac{7}{5}N_{j+1}^{s}}, & \mathcal{K}\cap \mathcal{Z}_{j}=\emptyset.
    \end{array}
  \right.
\end{equation}
(b) When $A_{j}\in\mathcal{RS}(40N^{s}N_{j+1},e^{-\frac{1}{50}N_{j}^{s}}) $ ,there is $\breve{k}_{j}\in\mathbb{Z}^{d}$ satisfying $|\breve{k}_{j}|\leq40N_{j+1}$, such that $|2\rho_{j}-\langle\breve{k}_{j},\alpha\rangle|\leq e^{-\frac{1}{50}N_{j}^{s}}$. Then, $P_{J}\in SU(1,1)$ with $P_{j}A_{j}P_{j}^{-1}=\left(
                                           \begin{array}{cc}
                                             e^{i\rho} & 0 \\
                                             0 & -e^{i\rho} \\
                                           \end{array}
                                         \right)$ and $B_{j}$ with $|B_{j}-I|_{h_{j}}\leq |\lambda|\|W_{E}\|e^{-\frac{3}{5}N_{j}^{s}}$ such that
\begin{equation}\label{4.32}
|Q_{\breve{k}_{j}}F_{j+1}Q_{-\breve{k}_{j}}|_{\tilde{h}_{j}}\leq2|P_{j}\tilde{F}_{j}P_{j}^{-1}|_{\tilde{h}_{j}},
\end{equation}
\begin{equation}\label{4.34}
  A_{j+1}\in\mathcal{NR}(40N^{s}N_{j+2},e^{-\frac{1}{50}N_{j+1}^{s}}).
\end{equation}
If we write
$$A_{j+1}=\left(
            \begin{array}{cc}
             * & b_{j+1} \\
              \bar{b_{j+1}} & * \\
            \end{array}
          \right),~~~P_{j}Ad(\tilde{B}_{j})\cdot(\lambda v_{j}W_{E})P_{j}^{-1}=\left(
                                         \begin{array}{cc}
                                         * & g_{j}^{+} \\
                                           \bar{g}_{j}^{+} & *\\
                                         \end{array}
                                       \right),$$
we have
\begin{equation}\label{4.37}
\left\{
    \begin{array}{ll}
       |b_{j+1}-\hat{g}_{j}^{+}(\breve{k}_{j})|\leq|\lambda|\|W_{E}\|e^{-\frac{3}{2}|k_{j}|^{s}}e^{-|\breve{k}_{j}|h_{j}}, & \mathcal{K}\cap \mathcal{Z}_{j}\neq\emptyset; \\
      |b_{j+1}|\leq|\lambda|\|W_{E}\| e^{-\frac{3}{2}N_{j+1}^{s}}e^{-\frac{3}{4}|\breve{k}_{j}|h_{j-1}}, & \mathcal{K}\cap \mathcal{Z}_{j}=\emptyset.
    \end{array}
  \right.
\end{equation}
\end{Lemma}

\begin{pf}

Assume that the first step $j$ iteration is valid ,that is,~$B_{j}\in \mathcal{B}_{h_{j}}^{\omega}(\mathbb{T}^{d},SL(2,\mathbb{R}))$ has been constructed to make
$$\widetilde{B}_{j}(\theta+\alpha)A_{0}e^{\widetilde{F}(\theta)}\widetilde{B}_{j}^{-1}(\theta)=A_{j}e^{F_{j}(\theta)} .$$
We will use Proposition \ref{lem3.1} to construct $B_{j}$ such that it conjugates $Cocycle~(\alpha,A_{j}e^{\tilde{F}_{j}})$ to $Cocycle~(\alpha,A_{j+1}e^{{F}_{j+1}}).$

Let
   \begin{equation}\label{4.5}
   \tilde{h}_{j}= \left\{
    \begin{array}{ll}
    h_{j},&K\cap Z_{j}\neq\emptyset;\\
    \frac{3}{4}h_{j-1},&K\cap Z_{j}=\emptyset
    \end{array}
    \right.
   \end{equation}
 There are two different cases of $\tilde{F}_{j}$.\\
(1)$\mathcal{K}\cap \mathcal{Z}_{j}\neq\emptyset:$~by \eqref{4.1} and \eqref{6},~we have
\begin{equation}\label{4.6}
\begin{array}{ll}
    |\lambda Ad(\widetilde{B}_{j}).(v_{j}W_{E})|_{\tilde{h}_{j}}&\leq |\lambda|\|W_{E}\| e^{-\frac{9}{10}|k_{j}|^{s}}e^{\frac{1}{20}N_{j}^{s}}\\
    &\leq|\lambda|\|W_{E}\| e^{-\frac{17}{20}|k_{j}|^{s}},
    \end{array}
\end{equation}
\begin{equation}\label{4.7}
\begin{array}{ll}
  |\tilde{F}_{j}|_{\tilde{h}_{j}}&\leq |\lambda|\|W_{E}\| e^{-N_{j+1}^{s}}+|\lambda|\|W_{E}\| e^{-\frac{17}{20}|k_{j}|^{s}} \\
  &\leq|\lambda|\|W_{E}\| e^{-\frac{4}{5}|k_{j}|^{s}},
\end{array}
\end{equation}
\begin{equation}\label{4.8}
\begin{array}{ll}
    \langle\tilde{F}_{j}\rangle&=\langle F_{j}+\lambda Ad(\widetilde{B}_{j}).(v_{j}W_{E})\rangle \\
    &=\langle \lambda Ad(\widetilde{B}_{j}).(v_{j}W_{E})\rangle
  \end{array}
\end{equation}
(2)$\mathcal{K}\cap \mathcal{Z}_{j}=\emptyset:$~by \eqref{5},~$v_{j}=0$ and $\tilde{F}_{j}=F_{j}$,~we have
\begin{equation}\label{4.9}
|\tilde{F}_{j}|_{\tilde{h}_{j}}=|F_{j}|_{\frac{3}{4}h_{j-1}}\leq |\lambda|\|W_{E}\| e^{-N_{j+1}^{s}},
\end{equation}
\begin{equation}\label{4.10}
\langle\tilde{F}_{j}\rangle= \langle F_{j}\rangle=0.
\end{equation}
In both cases,~we have by \eqref{6} that
\begin{equation}\label{4.11}
|\tilde{F}_{j}|_{\tilde{h}_{j}}\leq |\lambda|\|W_{E}\| e^{-\frac{4}{5}N_{j}^{s}}.
\end{equation}

(a) $A_{j}\in\mathcal{NR}(40N^{s}N_{j+1},e^{-\frac{1}{50}N_{j}^{s}}) $:~we use Proposition \ref{lem3.2} to construct $B_{j}=e^{Y_{j}}$ conjugating $Cocycle ~(\alpha,A_{j}e^{\tilde{F}_{j}})$  to $Cocycle ~(\alpha,A_{j+1}e^{{F}_{j+1}}),$ where $A_{j+1}\in SL(2,\mathbb{R}),Y_{j}\in \mathcal{B}_{h_{j}}(\mathbb{T}^{d},SL(2,\mathbb{R}))$ and $F_{j+1}\in \mathcal{B}_{\tilde{h}_{j}}^{\geq40N^{s}N_{j+1}}(\mathbb{T}^{d},SL(2,\mathbb{R}).$ Together with \eqref{4.11},~we have the following estimates
\begin{eqnarray*}
&&|Y_{j}|_{\tilde{h}_{j}}\leq10e^{\frac{4}{50}N_{j}^{s}}|\tilde{F}_{j}|_{\tilde{h}_{j}},\quad
|F_{j+1}|_{\tilde{h}_{j}}\leq2|\tilde{F}_{j}|_{\tilde{h}_{j}},\\
&&\|A_{j}e^{\langle\tilde{F}_{j}\rangle}-A_{j+1}\|\leq400e^{\frac{4}{50}N_{j}^{s}}|\tilde{F}_{j}|_{\tilde{h}_{j}}^{2}.
\end{eqnarray*}
By (\ref{4.7}-\ref{4.8}),~we have
\begin{eqnarray}\label{4.12}
&&|Y_{j}|_{h_{j}}\leq10|\lambda|\|W_{E}\| e^{\frac{4}{50}N_{j}^{s}}e^{-\frac{4}{5}|k_{j}|^{s}}\leq\frac{1}{2}|\lambda|\|W_{E}\| e^{-\frac{3}{5}N_{j}^{s}},\\
\label{4.13}
&&|F_{j+1}|_{h_{j}}\leq2|\tilde{F}_{j}|_{\tilde{h}_{j}}\leq2|\lambda|\|W_{E}\| e^{-\frac{4}{5}N_{j}^{s}}
\end{eqnarray}
When $\mathcal{K}\cap \mathcal{Z}_{j}\neq\emptyset,$~we have
$$\begin{array}{ll}
  \|A_{j}e^{\langle \lambda Ad(\widetilde{B}_{j}).(v_{j}W_{E})\rangle}-A_{j+1}\|&\leq400|\lambda|\|W_{E}\| e^{\frac{4}{50}N_{j}^{s}}e^{-\frac{8}{5}|k_{j}|^{s}}\\
  &\leq|\lambda|\|W_{E}\| e^{-\frac{7}{5}|k_{j}|^{s}}
  \end{array}$$ by  (\ref{4.7}-\ref{4.8}).
When $\mathcal{K}\cap \mathcal{Z}_{j}=\emptyset,$
$$\|A_{j}-A_{j+1}\|\leq400|\lambda|\|W_{E}\| e^{\frac{4}{50}N_{j}^{s}}e^{-2N_{j+1}^{s}}\leq\lambda|\|W_{E}\| e^{-\frac{7}{5}N_{j+1}^{s}}$$
by (\ref{4.9}-\ref{4.10})
Then
\begin{equation}\label{4.14b}
\left\{
    \begin{array}{ll}
      \|A_{j+1}-I\|\geq\|A_{j}e^{\langle \lambda Ad(\widetilde{B}_{j}).(v_{j}W_{E})\rangle}-I\|-|\lambda|\|W_{E}\| e^{-\frac{7}{5}|k_{j}|^{s}} , & \mathcal{K}\cap \mathcal{Z}_{j}\neq\emptyset; \\
      \|A_{j+1}-I\|\geq\|A_{j}-I\|-|\lambda|\|W_{E}\| e^{-\frac{7}{5}N_{j+1}^{s}}, & \mathcal{K}\cap \mathcal{Z}_{j}=\emptyset.
    \end{array}
  \right.
\end{equation}
 By \eqref{4.6}\eqref{4.12},~we have
\begin{equation}\label{4.15b}
|B_{j}-I|_{h_{j}}\leq|\lambda|\|W_{E}\| e^{-\frac{3}{5}N_{j}^{s}},\|A_{j}-A_{j+1}\|\leq |\lambda|\|W_{E}\| e^{-\frac{3}{5}N_{j}^{s}}
\end{equation}
and
\begin{equation}\label{4.16}
\begin{array}{ll}
   |F_{j+1}|_{\frac{3}{4}h_{j}}&=\sum_{|k|\geq40N^{s}N_{j+1}}\|\hat{F}_{j}(k)\|e^{-\frac{3}{4}|k|h_{j}}\\
   &\leq e^{-40N^{s}N_{j+1}\times\frac{1}{4}h_{j}}\sum_{|k|\geq40N^{s}N_{j+1}}\|\hat{F}_{j}(k)\|e^{-|k|h_{j}}\\
   &\leq e^{-10N^{s}N_{j+1}h_{j}}\times 2|\lambda|\|W_{E}\| e^{-\frac{4}{5}N_{j}^{s}}\\
   &\leq |\lambda|\|W_{E}\| e^{-N_{j+2}^{s}}
   \end{array}
\end{equation}
Then  \eqref{4.1}  follows from  (\ref{4.14}-\ref{4.16}) and $\langle F_{j+1}\rangle=0$

(b) $A_{j}\in\mathcal{RS}(40N^{s}N_{j+1},e^{-\frac{1}{50}N_{j}^{s}}) $:~In this case,~there is $\breve{k}_{j}\in\mathbb{Z}^{d}$ satisfying $|\breve{k}_{j}|\leq40N_{j+1}$,~such that $|2\rho_{j}-\langle\breve{k}_{j},\alpha\rangle|\leq e^{-\frac{1}{50}N_{j}^{s}}.$
Then,~by \eqref{9}
\begin{equation}\label{4.28}
|2\rho_{j}|\in[e^{-\frac{1}{2000}N^{s-1}N_{j}^{s}},20N^{s}N_{j+1}(|\alpha|+1)],
\end{equation}
one can find $P_{j}\in SU(1,1)$ satisfying
\begin{equation}\label{4.29}
(1+|\rho_{j}|)\|P_{j}\|^{2}\leq(2|\rho_{j}|+4|\rho_{j}|^{-1}+4)\leq e^{\frac{1}{1500}N^{s-1}N_{j}^{s}},
\end{equation}
such that $P_{j}A_{j}P_{j}^{-1}=\left(
                          \begin{array}{cc}
                            e^{i\rho} & 0 \\
                            0 & -e^{i\rho} \\
                          \end{array}
                        \right).~$
It is obvious that $P_{j}$ conjugates $cocycle~(\alpha,A_{j}e^{\tilde{F}_{j}})$ to $cocycle~(\alpha,A_{j}'e^{g_{j}})$ with
$$A_{j}'=\left(
                          \begin{array}{cc}
                            e^{i\rho} & 0 \\
                            0 & -e^{i\rho} \\
                          \end{array}
                        \right),\quad g_{j}={P_{j}\tilde{F}_{j}P_{j}^{-1}}.$$
And $g_{j}$ satisfies
\begin{equation}\label{4.30}
|g_{j}|_{\tilde{h}_{j}}\leq|P_{j}\tilde{F}_{j}P_{j}^{-1}|_{\tilde{h}_{j}}\leq|\lambda|\|W_{E}\| e^{-\frac{7}{10}N_{j}^{s}}.
\end{equation}
By Proposition \ref{lem3.1}, there exists $Y_{j}\in \mathcal{B}_{\tilde{h}_{j}}^{\alpha}(\mathbb{T}^{d},SL(2,\mathbb{R})$ conjugates $cocycle~(\alpha,A_{j}'e^{g_{j}})$ to $cocycle~(\alpha,A_{j}'e^{g_{j}^{*}})$ i.e.
$$e^{Y_{j}(\theta+\alpha)}(A_{j}'e^{g_{j}(\theta)})e^{-Y_{j}(\theta)}=A_{j}'e^{g_{j}^{*}(\theta)}$$
with
\begin{equation}\label{b1}
  |Y_{j}|_{\tilde{h}_{j}}\leq2e^{\frac{1}{50}N_{j}^{s}}|g_{j}|_{\tilde{h}_{j}},\quad |g_{j}^{*}-g_{j}^{(re)}|\leq2e^{\frac{7}{50}N_{j}^{s}}.
\end{equation}
By \eqref{4.30} and \eqref{b1},we have
\begin{equation}\label{b2}
  |g_{j}^{*}|_{\tilde{h}_{j}}\leq2|g_{j}|_{\tilde{h}_{j}}\leq2|\lambda|\|W_{E}\| e^{-\frac{7}{10}N_{j}^{s}},
\end{equation}
with
\begin{equation}\label{b3}
  g_{j}^{*}(\theta)=\left(
                      \begin{array}{cc}
                        ia_{j} & b_{j}^{*}e^{i\langle \breve{k}_{j},\theta\rangle} \\
                        \bar{b}_{j}^{*}e^{i\langle \breve{k}_{j},\theta\rangle} &  -ia_{j} \\
                      \end{array}
                    \right)+\mathcal{R}_{N}g_{j}^{*},
\end{equation}
where $a_{j}\in\mathbb{R}$,$b_{j}^{*}\in\mathbb{C}$ with the estimates
\begin{equation}\label{b4}
  |a_{j}|\leq|g_{j}^{*}|_{\tilde{h}_{j}}\leq2|\lambda|\|W_{E}\| e^{-\frac{7}{10}N_{j}^{s}},
  \quad |b_{j}^{*}|\leq|g_{j}^{*}|_{\tilde{h}_{j}}e^{-|\breve{k}_{j}|\tilde{h}_{j}}\ll2|\lambda|\|W_{E}\| e^{-\frac{7}{10}N_{j}^{s}}.
\end{equation}

Using Proposition  \ref{lem3.2} (we choose $\tilde{N}=200N^{s}N_{j+1},~N=40N^{s}N_{j+1}$),~together with \eqref{4.30} and \eqref{4.11},~we construct $A_{j+1}\in SL(2,\mathbb{R}),~Y_{j}\in \mathcal{B}_{\tilde{h}_{j}}^{\alpha}(\mathbb{T}^{d},SL(2,\mathbb{R})$ and $F_{j+1}\in \mathcal{B}_{\tilde{h}_{j}}^{\geq40N^{s}N_{j+1}}(\mathbb{T}^{d},SL(2,\mathbb{R}),$~such that $Q_{-\breve{k}_{j}}e^{Y_{j}}$ conjugates \\ $cocycle~(\alpha,\left(
                          \begin{array}{cc}
                            e^{i\rho} & 0 \\
                            0 & -e^{i\rho} \\
                          \end{array}
                        \right)e^{P_{j}\tilde{F}_{j}P_{j}^{-1}})$ to $cocycle~(\alpha,A_{j+1}e^{F_{j+1}})$ ,where $A_{j+1}=\tilde{A}_{j}^{'}e^{H_{j}}=e^{A_{j+1}^{''}}$(using \eqref{3.6}).
Then,~by (\ref{4.29}-\ref{4.31}) we get
\begin{equation}\label{4.33}
|Q_{\breve{k}_{j}}B_{j}P_{j}^{-1}-I|_{h_{j}}\leq|\lambda|\|W_{E}\|e^{-\frac{3}{5}N_{j}^{s}},\quad \|P_{j}\|^{2}\leq\frac{1}{1+\rho_{j}}e^{\frac{1}{1500}N^{s-1}N_{j}^{s}},
\end{equation}
and
\begin{equation}\label{4.31}
\|A_{j+1}^{''}\|\leq e^{-\frac{1}{50}N_{j}^{s}},\quad |Y_{j}|_{\tilde{h}_{j}}\leq10e^{\frac{1}{50}N_{j}^{s}}|P_{j}\tilde{F}_{j}P_{j}^{-1}|_{\tilde{h}_{j}}
\end{equation}
If denote $spec(A_{j+1}^{''})=\pm\mu_{j+1}$ and
$A_{j+1}^{''}\triangleq\left(
                         \begin{array}{cc}
                           ic_{j+1} & d_{j+1} \\
                           \bar{d}_{j+1} & -ic_{j+1} \\
                         \end{array}
                       \right)$, then
$$|\mu_{j+1}|\leq\sqrt{||c_{j+1}|^{2}-|d_{j+1}|^{2}|}\leq\sqrt{2}e^{-\frac{1}{50}N_{j}^{s}},$$
which implies that
$$|\rho_{j+1}|\triangleq|rot(\alpha,A_{j+1})|\leq2e^{-\frac{1}{50}N_{j}^{s}},$$
so
$$|2\rho_{j+1}-\langle k,\alpha\rangle|\geq|\langle k,\alpha\rangle|-2e^{-\frac{1}{50}N_{j}^{s}}\geq e^{-\frac{1}{50}N_{j}^{s}},\quad
0<|k|\leq40N^{s}N_{j+2},$$
then
\begin{equation}\label{4.34b}
  A_{j+1}\in\mathcal{NR}(40N^{s}N_{j+2},e^{-\frac{1}{50}N_{j+1}^{s}}).
\end{equation}

Denote $spec(A_{j+1}^{''})=\{\pm\mu_{j}\}$, by \eqref{b4}, we have
\begin{equation}\label{b5}
  \mu_{j}^{2}=|a_{j}|^{2}-|b_{j}^{*}|^{2}\leq|g_{j}^{*}|_{\tilde{h}_{j}}.
\end{equation}
For $A_{j+1}=\tilde{A}_{j}e^{H_{j}}$, then by \eqref{A2-trans} in Lemma \ref{Lem-A1} we get
$$b_{j+1}=b_{j}^{*}\mu_{j}^{-1}sinh(\mu_{j})e^{i(\rho_{j}-\frac{\langle \breve{k}_{j},\alpha\rangle}{2})}.$$
When $\mathcal{K}\cap \mathcal{Z}_{j}=\emptyset$,~we know $P_{j}\tilde{F}_{j}P_{j}^{-1}=P_{j}F_{j}P_{j}^{-1}$. For all $j\geq1$ ,recall $|2\rho_{j+1}-\langle \breve{k}_{j},\alpha\rangle|<e^{-\frac{1}{50}N_{j}^{s}},~~g_{j}=P_{j}\tilde{F}_{j}P_{j}^{-1}$ and by \eqref{11} we get
\begin{equation}\label{4.35}
  |b_{j}^{*}|\leq|\lambda|\|W_{E}\| e^{\frac{1}{1500}N^{s-1}N_{j}^{s}-N_{j+1}^{s}}e^{-\frac{3}{4}|\breve{k}_{j}|h_{j-1}}.
\end{equation}
By \eqref{b2},\eqref{b4}, \eqref{b5}and \eqref{4.35} we get
\begin{equation}\label{4.36}
\begin{array}{ll}
    |b_{j}-b_{j}^{*}|&=b_{j}^{*}(\mu_{j}^{-1}sinh(\mu_{j})e^{i(\rho_{j}-\frac{\langle \breve{k}_{j},\alpha\rangle}{2})}-1)\\
    &\leq |g_{j}^{*}|^{2}_{\tilde{h}_{j}}e^{-|\breve{k}_{j}|h_{j-1}}\\
    &\leq 400e^{\frac{1}{50}N_{j}^{s}+\frac{2}{1500}N^{s-1}N_{j}^{s}}e^{-2N_{j+1}^{s}}e^{-\frac{3}{4}|\breve{k}_{j}|h_{j-1}}\\
    &<\frac{1}{2}|\lambda|\|W_{E}\| e^{-\frac{3}{4}|\breve{k}_{j}|h_{j-1}},
  \end{array}
\end{equation} for the inequalities $\frac{sinh(x)}{x}\leq1+x^{2}, e^{x}\leq1+2x.$
Then, by \eqref{b1} have
\begin{equation}\label{b7}
  |\hat{g}_{j}^{*}(\breve{k}_{j})-\hat{g}_{j}(\breve{k}_{j})|_{\tilde{h}_{j}}\leq2e^{\frac{7}{50}N_{j}^{s}}|\tilde{F}_{j}|^{2}_{\tilde{h}_{j}}
  \|P_{j}\|^{4}e^{-|\breve{k}_{j}|\tilde{h}_{j}},
\end{equation}
and
\begin{equation}\label{b8}
  |g_{j}(x)-\lambda P_{j}Ad(\widetilde{B}_{j})v_{j}(\theta)WP_{j}^{-1}|_{\tilde{h}_{j}}\leq|\lambda|\|W_{E}\||\tilde{F}_{j}|_{\tilde{h}_{j}}
  \|P_{j}\|^{2}\ll\frac{1}{2}|\lambda|\|W_{E}\|e^{-\frac{3}{2}N_{j+1}^{s}}.
\end{equation}
Together \eqref{4.1}, \eqref{4.30},\eqref{b2}, and \eqref{4.36}-\eqref{b8}, get \eqref{4.37}.\\

Now we verify (\ref{4.1}).~By \eqref{4.32}
\begin{equation}\label{4.38}
\begin{array}{ll}
       |F_{j+1}|_{\frac{3}{4}h_{j}}&\leq|Q_{\breve{k}_{j}}|_{\frac{3}{4}h_{j}}^{2}\|F_{j,+}\|_{\frac{3}{4}h_{j}} \\
       &\leq|Q_{\breve{k}_{j}}|_{\frac{3}{4}h_{j}}^{2}\sum_{|k|\geq200N^{s}N_{j+1}}\|\hat{F}_{j,+}(k)\|e^{\frac{3}{4}|k|_{h_{j}}} \\
       &\leq e^{-20N^{s}N_{j+1}h_{j}}\sum_{|k|\geq200N^{s}N_{j+1}}\|\hat{F}_{j,+}(k)\|e^{|k|_{h_{j}}} \\
       &\leq e^{-20N^{s}N_{j+1}h_{j}} \times2|\lambda|\|W_{E}\| e^{-\frac{7}{10}N_{j}^{s}}\\
       &\leq|\lambda|\|W_{E}\| e^{-N_{j+2}^{s}},
     \end{array}
\end{equation}
where $F_{j,+}=Q_{\breve{k}_{j}}F_{j+1}Q_{\breve{k}_{j}}^{-1}$.~Note that $Q_{\breve{k}_{j}}$ satisfies
$$|Q_{\breve{k}_{j}}|_{h_{j+1}}\leq e^{\frac{1}{2}\times40N^{s}N_{j+1}h_{j+1}}\leq e^{2N^{2s-1}N_{j+1}^{s}}\leq e^{\frac{1}{100}N_{j+1}^{s}},$$
and then $B_{j}=Q_{-\breve{k}_{j}}e^{Y_{j}}P_{j}$ has the estimation
$$|B_{j}|_{h_{j+1}}\leq2e^{\frac{1}{3000}N_{j}^{s}}e^{\frac{1}{100}N_{j+1}^{s}}\leq e^{\frac{1}{80}N_{j+1}^{s}}.$$
By\eqref{4.1},~we have
\begin{equation}\label{4.39}
  \begin{array}{ll}
    |\tilde{B}_{j+1}|_{h_{j+1}}&\leq|\tilde{B}_{j}|_{h_{j}}|B_{j}|_{h_{j+1}} \\
    &\leq e^{\frac{1}{40}N_{j}^{s}}e^{\frac{1}{80}N_{j+1}^{s}} \leq e^{\frac{1}{40}N_{j+1}^{s}}.
  \end{array}
\end{equation}
By \eqref{4.34},(\ref{4.38}-\ref{4.39}) and the fact $\langle F_{j+1}\rangle=0$,~we get \eqref{4.1} .
\end{pf}

\section{PROOF OF  THEOREM \ref{th-1.2}}
In section 3 it is stated  that $cocycle ~(\alpha,Ae^{\Sigma_{j=1}^{\infty}\lambda v_{j}(\theta)W})$ has been conjugated to cocycle $$(\alpha,A_{j+1}e^{F_{j+1}+ Ad(\tilde{B}_{j+1}).(\Sigma_{p=j+1}^{\infty}\lambda v_{p}(\theta)W_{E})})$$ by $\tilde{B}_{j+1}$.~We now assume that the fibered rotation number $\rho(E)$ of the  $Cocyle~(\alpha,Ae^{\Sigma_{j+1}^{\infty}\lambda v_{j}(\theta)W})$ is $\frac{1}{2}\langle k_{J},\omega\rangle$ for arbitrary and fixed $k_{J}\in\mathcal{K}$.~The idea is to prove that,~in Lemma \ref{lem4.1},~the fibered rotation number of $A_{j}$ is zero and $\|A_{j}-I\|$ is uniformly bounded away from zero for sufficiently large $j$.~This implies that the gap with the labelling $k_{J}$ is open.

\subsection{More Estimates}
Let us introduce some notation.~Define
\begin{equation}\label{5.1}
W_{j}\triangleq Ad(B_{j}).W_{E}=B_{j}W_{E}(B_{j})^{-1},
\end{equation}
Let us define the integer $j_{0}$ as
\begin{equation}\label{5.2}
  j_{0}\triangleq min\{1\leq j\leq J|A_{j}\in\mathcal{RS}(40N^{s}N_{j+1},e^{-\frac{1}{50}N_{j}^{s}})\}
\end{equation}
And for $j_{0}<j<J,$~we write $W_{j}$ as
$$W_{j}=\left(
          \begin{array}{cc}
            iu_{j} & w_{j}e^{-i\langle\widetilde{{k}_{j}},\cdot\rangle} \\
          \overline{w_{j}}e^{-i\langle\widetilde{{k}_{j}},\cdot\rangle} & -iu_{j} \\
          \end{array}
        \right),
$$
and define quantities
\begin{equation}\label{5.3}
\xi_{j}\triangleq|\langle w_{j}\rangle|,~~\mathcal{M}_{j}\triangleq|w_{j}|_{h_{j}}+|u_{j}|_{h_{j}},
\end{equation}
\begin{equation}\label{5.4}
\mathfrak{m}_{j}\triangleq\sup_{k\in\mathbb{Z}^{d},|k|>|\widetilde{{k}_{j}}|}\frac{1}{2}(|\widehat{w_{j}}(k)|+|\widehat{u_{j}}(k)|)
\end{equation}
The lower bounds estimation of $\|A_{j}-I\|$  will play a crucial  role in the proof of Theorem \ref{th-1.2}.~And the quantity $\xi_{j}$ for $j\geq j_{0}$ are indispensable.~However,~in every step of KAM, the quantities $\mathcal{M}_{j},~\mathfrak{m}_{j}$ and $\xi_{j}$ will influence the estimate of $\xi_{j}$ .~Thus,~we need to estimate $\xi_{j},~\mathcal{M}_{j},~\mathfrak{m}_{j}$ for $j>j_{0}$ in every step.

Next we will estimate $\breve{k}_{j}$ and $\tilde{k}_{j}$.

\begin{Lemma}\cite[Lemma 6.1]{HSY}\label{lem5.1}
$\forall$ $j\geq1$,~we have
\begin{equation}\label{5.5}
|\tilde{k}_{j}|\leq41N^{s}N_{j+1}.
\end{equation}
For $j>1$ and $A_{j}\in\mathcal{RS}(40N^{s}N_{j+1},e^{-\frac{1}{50}N_{j}^{s}})$,~we have
\begin{equation}\label{5.6}
|\tilde{k}_{j}|\leq41N^{s}N_{j-1},
\end{equation}
\begin{equation}\label{5.7}
|\breve{k}_{j}|>40N^{s}N_{j},
\end{equation}
\begin{equation}\label{5.8}
|\tilde{k}_{j+1}|>39N^{s}N_{j}.
\end{equation}
As a consequence,~for $j>1$ and $A_{j}\in\mathcal{RS}(40N^{s}N_{j+1},e^{-\frac{1}{50}N_{j}^{s}})$
\begin{equation}\label{5.9}
\begin{array}{ll}
&e^{|\tilde{k}_{j}|h_{j}}\leq e^{5N^{2s-2}N_{j}^{s}},\\
&e^{-|\breve{k}_{j}|h_{j}}\leq e^{-4N^{2s-1}N_{j}^{s}}\\
&e^{-|\tilde{k}_{j+1}|h_{j}}\leq e^{-\frac{399}{100}N^{2s-1}N_{j}^{s}}.
\end{array}
\end{equation}
\end{Lemma}

\begin{pf}
We prove \eqref{5.5} inductively.~When $j=1$ it is obviously true  .~Assume that \eqref{5.5} is true for the step $j$.~If $A_{j}\in\mathcal{NR}(40N^{s}N_{j+1},e^{-\frac{1}{50}N_{j}^{s}})$,~\eqref{5.5} with $j+1$ is obvious since $\tilde{k}_{j+1}=\tilde{k}_{j}$.~If
$A_{j}\in\mathcal{RS}(40N^{s}N_{j+1},e^{-\frac{1}{50}N_{j}^{s}})$,~we have

$$
\begin{array}{ll}
  |\tilde{k}_{j+1}|&\leq|\tilde{k}_{j}|+|\breve{k}_{j}| \\
  &\leq 100N^{s}N_{j}+40N^{s}N_{j+1}\\
  &\leq 100N^{s}N_{j+1},
\end{array}
$$
which also verifies \eqref{5.5}.
If $j>1$ and $A_{j}\in\mathcal{RS}(40N^{s}N_{j+1},e^{-\frac{1}{50}N_{j}^{s}})$,~in Proposition\ref{lem4.1},~$A_{j-1}\in\mathcal{NR}(40N^{s}N_{j},e^{-\frac{1}{50}N_{j-1}^{s}})$.~Then,~by \eqref{5.5} and $\tilde{k}_{j}=\tilde{k}_{j-1}$,~we get \eqref{5.6}.~We will prove \eqref{5.7} via contradiction.~If $\breve{k}_{j}\leq40N^{s}N_{j}$,~by \eqref{4.15},
 $$ |\rho_{j}-\rho_{j-1}|\leq15(\|A_{j}-A_{j-1}\|)^{\frac{1}{2}}\leq15e^{-\frac{3}{10}N_{j-1}^{s}},$$
$$
  |2\rho_{j}-\langle\breve{k}_{j},\omega\rangle|<e^{-\frac{1}{50}N_{j}^{s}},
$$
and then,~it follows that
$$
|2\rho_{j-1}-\langle\breve{k}_{j},\omega\rangle|<e^{-\frac{1}{50}N_{j}^{s}}+30e^{-\frac{3}{10}N_{j-1}^{s}}<e^{-\frac{1}{50}N_{j-1}^{s}},
$$
which contradicts $A_{j-1}\in\mathcal{NR}(40N^{s}N_{j},e^{-\frac{1}{50}N_{j-1}^{s}})$.~For $\tilde{k}_{j+1}=\tilde{k}_{j}+\breve{k}_{j}$,~by (\ref{5.6}-\ref{5.7}), ~$|\tilde{k}_{j}|<\frac{1}{1000}|\breve{k}_{j}|$,~and then,~one can get \eqref{5.8}.~By (\ref{5.6}-\ref{5.8}) and \eqref{2} (the definition of $h_{j}$),~\eqref{5.9} follows.
\end{pf}

 Then from the definition of $j_{0}$ in \eqref{5.2},~it is obvious that $\tilde{k}_{j}=0$ for $j<j_{0}$ and $B_{j}$ is close to the identity for $1\leq j<j_{0}$ as $j_{0}>1$.~It follows that
\begin{Lemma}\cite[Lemma 6.2]{HSY}\label{lem5.2}
For $1\leq j\leq j_{0}$
\begin{equation}\label{5.10}
\begin{array}{ll}
 & \|A_{j}\|\leq2 \|W_{E}\|^{-1},\\
 &|W_{j}|_{h_{j}}\leq 2\|W_{E}\|,\\
 &\|[A_{j},\langle W_{j}\rangle]\|\geq\frac{1}{2}.
  \end{array}
\end{equation}
\end{Lemma}

\begin{pf}
Note that,~by Lemma \ref{lem4.1},~for $j<j_{0}$
$$
  \|A_{j+1}-A_{j}\|\leq\|W_{E}\| e^{-\frac{3}{5}N_{j}^{s}},\quad |W_{j+1}-W_{j}|_{h_{j}}\leq4\|W_{E}\| e^{-\frac{3}{5}N_{j}^{s}}
$$
where $W_{j+1}=Ad(\tilde{B}_{j}).W_{j}$ with $|\tilde{B}_{j}-I|_{h_{j}}\leq\|W_{E}\| e^{-\frac{3}{5}N_{j}^{s}}$.
\end{pf}

$\mathcal{M}_{j},~\mathfrak{m}_{j},~\xi_{j}$ defined in (\ref{5.3}- \ref{5.4}) have the following estimates

\begin{Lemma}\cite[Lemma 6.3]{HSY}\label{lem5.3}
\begin{equation}\label{5.11}
\begin{array}{ll}
  &\mathcal{M}_{j_{0}+1}\leq e^{\frac{1}{1000}N^{s-1}N_{j_{0}}^{s}},\\
  &\mathfrak{m}_{j_{0}+1}\leq e^{-\frac{399}{100}N^{2s-1}N_{j_{0}}^{s}},\\
  &\xi_{j_{0}+1}\geq e^{-\frac{1}{1000}N^{s-1}N_{j_{0}}^{s}}.
\end{array}
\end{equation}
\end{Lemma}

\begin{pf}
Notice that $A_{j_{0}}\in\mathcal{RS}(40N^{s}N_{j_{0}+1},e^{-\frac{1}{50}N_{j_{0}}^{s}})$ and $B_{j_{0}}=Q_{-\breve{k}_{j_{0}}}e^{Y_{j_{0}}}P_{j_{0}}^{-1}$ in Lemma \ref{lem4.1} have estimates
\begin{equation}\label{5.12}
  |e^{Y_{j_{0}}}-I|_{h_{j_{0}}}\leq|\lambda| e^{-\frac{3}{5}N_{j_{0}}^{s}},\quad \|P_{j_{0}}\|^{2}\leq\frac{1}{1+|\rho_{j_{0}}|}e^{\frac{1}{1500}}N^{s-1}N_{j_{0}}^{s}.
\end{equation}
Write
$$
  Ad(e^{Y_{j_{0}}}P_{j_{0}}^{-1}).W_{j_{0}}=\left(
                                              \begin{array}{cc}
                                                iu_{j_{0}+1} & w_{j_{0}+1} \\
                                                 \overline{w_{j_{0}+1}} & -iu_{j_{0}+1} \\
                                              \end{array}
                                            \right),
$$
$$
  W_{j_{0}+1}=Ad(Q_{-\breve{k}_{j_{0}}}).\left(
                                              \begin{array}{cc}
                                                iu_{j_{0}+1} & w_{j_{0}+1} \\
                                                 \overline{w_{j_{0}+1}} & -iu_{j_{0}+1} \\
                                              \end{array}
                                            \right).
$$
Then,~we get $(|W|_{h_{j_{0}}})\leq2$
$$
  \mathcal{M}_{j_{0}+1}\leq4\|P_{j_{0}}\|^{2}(1+4|e^{Y_{j_{0}}}-I|_{h_{j_{0}}})\leq20\|P_{j_{0}}\|^{2}\leq e^{\frac{1}{1000}}N^{s-1}N_{j_{0}}^{s}.
$$
By Lemma \ref{lem5.2},~$\|[A_{j_{0},\langle W_{j_{0}}\rangle}]\|\geq\frac{1}{2}$.~Then,~using Lemmas \ref{lem2.5},\ref{lem2.6} and \eqref{5.12},~we get
$$
\begin{array}{ll}
  \xi_{j_{0}+1}&\geq\frac{1}{4\rho_{j_{0}}\|P_{j_{0}}\|^{2}}-4|e^{Y_{j_{0}}}-I|_{h_{j_{0}}}\|P_{j_{0}}\|^{2}|W|_{h_{j_{0}}}\\
  &\geq e^{-\frac{1}{1000}}N^{s-1}N_{j_{0}}^{s}.
\end{array}
$$
By Lemma \ref{lem4.1},~the desired estimations for $\mathcal{M}_{j_{0}+1}$ and $\xi_{j_{0}+1}$ hold.

To estimate $\mathfrak{m}_{j_{0}+1}$,~we consider two different cases.~When $j_{0}>1$,~
$$
\begin{array}{ll}
  \mathfrak{m}_{j_{0}+1}&\leq\mathcal{M}_{j_{0}+1}e^{-|\breve{k}_{j_{0}}|h_{j_{0}}}\\
  &\leq e^{\frac{1}{1000}}N^{s-1}N_{j_{0}}^{s}e^{-8N^{2s-1}N_{j_{0}^{s}}}\\
  &\leq e^{-7N^{2s-1}N_{j_{0}^{s}}}
  \end{array}
$$
by \ref{5.9}.~When $j_{0}=1$,~by \eqref{5.12} and \eqref{11},
$$
  \begin{array}{ll}
    \mathfrak{m}_{2}&\leq|W_{+}-\langle W_{+}\rangle|_{h_{1}} \\
    &\leq4|e^{Y_{1}}-I|_{h_{1}}\|P_{1}\|^{2}\\
    &\leq e^{-\frac{1}{5}N_{1}^{s}}\\
    &\leq e^{-7N^{2s-1}N_{1}^{s}},
  \end{array}
$$
where $W_{+}\triangleq Ad(e^{Y_{1}}P_{1}^{-1}).W=\left(
                                                   \begin{array}{cc}
                                                     iu_{2} & w_{2} \\
                                                     \overline{w_{2}} & -iu_{2}\\
                                                   \end{array}
                                                 \right)$.
\end{pf}

\begin{Lemma}\cite[Lemma 6.4]{HSY}\label{lem5.4}
Let $j_{0}+1\leq j\leq J$.~If $A_{j}\in\mathcal{NR}(40N^{s}N_{j+1},e^{-\frac{1}{50}N_{j}^{s}})$,~then
\begin{equation}\label{5.13}
\begin{array}{ll}
 & \mathcal{M}_{j+1}\leq\mathcal{M}_{j}(1+e^{-\frac{1}{10}N_{j}^{s}}),\\
 &\mathfrak{m}_{j+1}\leq \mathfrak{m}{j}+\mathcal{M}_{j}e^{-\frac{1}{10}N_{j}^{s}},\\
&\xi_{j+1}\geq\xi_{j}-\mathcal{M}_{j}e^{-\frac{1}{10}N_{j}^{s}}.
\end{array}
\end{equation}
And if $A_{j}\in\mathcal{RS}(40N^{s}N_{j+1},e^{-\frac{1}{50}N_{j}^{s}})$,~then
\begin{equation}\label{m-5.14}
\begin{array}{ll}
  &\mathcal{M}_{j+1}\leq\mathcal{M}_{j}e^{\frac{1}{1000}N^{s-1}N_{j}^{s}},\\
  &\mathfrak{m}_{j+1}\leq \mathcal{M}_{j}e^{-\frac{398}{100}N^{2s-1}N_{j}^{s}},\\
&\xi_{j+1}\geq\xi_{j}-3\mathfrak{m}{j}-\mathcal{M}_{j}e^{-\frac{1}{10}N_{j}^{s}}.
\end{array}
\end{equation}
\end{Lemma}

\begin{pf}
If $A_{j}\in\mathcal{NR}(40N^{s}N_{j+1},e^{-\frac{1}{50}N_{j}^{s}})$,~by Lemmas \ref{lem4.1} and \ref{lem5.1},
$$
|e^{Y_{j}}-I|_{h_{j}}\leq\sigma e^{-\frac{3}{5}N_{j}^{s}},\quad e^{|\tilde{k}_{j}|h_{j}}\leq e^{5N^{2s-1}N_{j}^{s}}.
$$
It follows that
$$
  \begin{array}{ll}
    |W_{j+1}-W_{j}|_{h_{j}}&\leq4|B_{j}-I|_{h_{j}}|W_{j}|_{h_{j}}\leq 4|B_{j}-I|_{h_{j}}e^{|\tilde{k}_{j}|h_{j}}\mathcal{M}_{j}\\
    &\leq4 e^{5N^{2s-1}N_{j}^{s}}e^{-\frac{3}{5}N_{j}^{s}}\mathcal{M}_{j}\leq e^{-\frac{1}{10}N_{j}^{s}}\mathcal{M}_{j}.
  \end{array}
$$
Equation \eqref{5.13} then follows from definition (\ref{5.3}-\ref{5.4}) and Lemma \ref{lem2.6}.
If $A_{j}\in\mathcal{RS}(40N^{s}N_{j+1},e^{-\frac{1}{50}N_{j}^{s}})$,~$B_{j}=Q_{-\breve{k}_{j}}\breve{B}_{j}P_{j}^{-1}$ ($\breve{B}_{j}=e^{Y_{j}}$)in Lemma \ref{lem4.1} satisfies
$$
  |\breve{B}_{j}-I|_{h_{1}}\leq|\lambda| e^{-\frac{3}{5}N_{j}^{s}},\quad  \|P_{j}\|^{2}\leq\frac{1}{1+|\rho_{j}|}e^{\frac{1}{1500}}N^{s-1}N_{j}^{s}.
$$
Same as the proof of Lemma \ref{lem5.3},~we get by using Lemma \ref{lem5.1} that
$$
  \begin{array}{ll}
  \mathcal{M}_{j+1}&\leq2 \|P_{j}\|^{2}(1+4|Q_{\breve{k}_{j}}\breve{B}_{j}P_{j}-I|_{h_{j}})|W_{j}|_{h_{j}}\\
    &\leq2e^{\frac{1}{1500}}N^{s-1}N_{j}^{s}(1+4|Q_{\breve{k}_{j}}\breve{B}_{j}P_{j}-I|_{h_{j}})e^{|\tilde{k}_{j}|h_{j}}\mathcal{M}_{j}\\
    &\leq4e^{\frac{1}{1500}}N^{s-1}N_{j}^{s}(1+4e^{-\frac{3}{5}N_{j}^{s}})e^{10N^{2s-2}N_{j}^{s}}\mathcal{M}_{j}\\
    &\leq e^{\frac{1}{1000}}N^{s-1}N_{j}^{s}\mathcal{M}_{j}.
  \end{array}
$$
By Lemmas \ref{lem2.5} and \ref{lem2.6},~we have
$$
  \begin{array}{ll}
    \xi_{j+1}&\geq\xi_{j}-3\mathfrak{m}_{j}-4\mathcal{M}_{j}|Q_{\breve{k}_{j}}\breve{B}_{j}P_{j}-I|_{h_{j}} \|P_{j}\|^{2}e^{|\tilde{k}_{j}|h_{j}} \\
    &\geq \xi_{j}-3\mathfrak{m}_{j}-4\mathcal{M}_{j}e^{-\frac{3}{5}N_{j}^{s}}e^{\frac{1}{1500}}N^{s-1}N_{j}^{s}e^{10N^{2s-2}N_{j}^{s}}\\
    &\geq \xi_{j}-3\mathfrak{m}_{j}-\mathcal{M}_{j}e^{-\frac{1}{10}N_{j}^{s}}.\\
    \mathfrak{m}_{j+1}&\leq \mathcal{M}_{j+1}e^{-|\tilde{k}_{j+1}|h_{j}}\leq \mathcal{M}_{j}e^{\frac{1}{1000}}N^{s-1}N_{j}^{s}e^{-7N^{2s-1}N_{j}^{s}}\\
    &\leq\mathcal{M}_{j}e^{-5N^{2s-1}N_{j}^{s}}.
  \end{array}
$$
\end{pf}

We now arrive at an important conclusion:
\begin{Lemma}\cite[Lemma 6.5]{HSY}\label{lem5.5}
For all $j>j_{0}$
\begin{equation}\label{m-5.15}
  \xi_{j}\geq10\mathfrak{m}_{j}+\frac{1}{2^{j}}\xi_{j_{0}+1},~~~for~ all ~j>j_{0}.
\end{equation}
\end{Lemma}

\begin{pf}
We first prove inductively the following
\begin{equation}\label{m-5.16}
  \mathcal{M}_{j}\leq e^{\frac{1}{500}N^{s-1}N_{j}^{s}(1-\frac{1}{2^{j}})},\quad for~ all j>j_{0}.
\end{equation}
By Lemma \ref{lem5.3},~$\mathcal{M}_{j_{0}+1}\leq e^{\frac{1}{1000}N^{s-1}N_{j_{0}}^{s}}$.~Assume that \eqref{m-5.16} holds for the step $j$.~If $A_{j}\in\mathcal{NR}(40N^{s}N_{j+1},e^{-\frac{1}{50}N_{j}^{s}})$,~by Lemma \ref{lem5.4} and \eqref{11}
$$
  \mathcal{M}_{j+1}\leq2e^{\frac{1}{500}N^{s-1}N_{j}^{s}(1-\frac{1}{2^{j}})}\leq e^{\frac{1}{500}N^{s-1}N_{j+1}^{s}(1-\frac{1}{2^{j+1}})}.
$$
If $A_{j}\in\mathcal{RS}(40N^{s}N_{j+1},e^{-\frac{1}{50}N_{j}^{s}})$,~by Lemma \ref{lem5.4} and \eqref{11}
$$
\begin{array}{ll}
\mathcal{M}_{j+1}&\leq e^{\frac{1}{500}N^{s-1}N_{j}^{s}(1-\frac{1}{2^{j}})}e^{\frac{1}{1000}N^{s-1}N_{j}^{s}}\leq e^{\frac{3}{1000}N^{s-1}N_{j}^{s}}\\
&e^{\frac{1}{500}N^{s-1}N_{j+1}^{s}(1-\frac{1}{2^{j+1}})}.
\end{array}
$$
We now prove \eqref{m-5.15} inductively.~By Lemma \ref{lem5.3},~\eqref{m-5.15} holds for $j_{0}+1$.~Assume that it holds for some $j\geq j_{0}+1$.~If $A_{j}\in\mathcal{NR}(40N^{s}N_{j+1},e^{-\frac{1}{50}N_{j}^{s}})$,~we have by Lemma \ref{lem5.4} and \eqref{m-5.16},
$$
  \xi_{j+1} \geq \xi_{j}-\mathcal{M}_{j}e^{-\frac{1}{10}N_{j}^{s}}\geq \xi_{j}-e^{-\frac{1}{20}N_{j}^{s}},$$
 $$ \mathfrak{m}_{j+1} \leq \mathfrak{m}_{j}+\mathcal{M}_{j}e^{-\frac{1}{10}N_{j}^{s}}\leq\mathfrak{m}_{j}+e^{-\frac{1}{20}N_{j}^{s}}.$$
With (by Lemma \ref{lem5.3} and \eqref{m-5.15}).
$$
  \begin{array}{ll}
    e^{-\frac{1}{20}N_{j}^{s}}&=e^{-\frac{1}{40}N_{j}^{s}}\times e^{-\frac{1}{40}N^{s}N_{j-1}^{s}} \\
    &\ll\frac{1}{11\times2^{j+1}}e^{-\frac{1}{1000}N^{s-1}N_{j_{0}}^{s}}\\
    &\leq\frac{1}{11\times2^{j+1}}\xi_{j_{0}+1},
  \end{array}
$$
we have
$$
  \begin{array}{ll}
    \xi_{j+1}&\geq\xi_{j}-e^{-\frac{1}{20}N_{j}^{s}}\geq10\mathfrak{m}_{j}+\frac{1}{2^{j}}\xi_{j_{0}+1}-e^{-\frac{1}{20}N_{j}^{s}} \\
    &\geq 10(\mathfrak{m}_{j+1}-e^{-\frac{1}{20}N_{j}^{s}})+\frac{1}{2^{j}}\xi_{j_{0}+1}-e^{-\frac{1}{20}N_{j}^{s}}\\
    &\geq 10\mathfrak{m}_{j+1}+\frac{1}{2^{j+1}}\xi_{j_{0}+1}+(\frac{1}{2^{j+1}}\xi_{j_{0}+1}-11e^{-\frac{1}{20}N_{j}^{s}})\\
    &\geq 10\mathfrak{m}_{j+1}+\frac{1}{2^{j+1}}\xi_{j_{0}+1}.
  \end{array}
$$
If $A_{j}\in\mathcal{RS}(40N^{s}N_{j+1},e^{-\frac{1}{50}N_{j}^{s}})$,~we get by Lemma \ref{lem5.4} and \eqref{m-5.16}
\begin{equation}\label{m-5.17}
  \xi_{j+1} \geq \xi_{j}-3\mathfrak{m}_{j}-\mathcal{M}_{j}e^{-\frac{1}{10}N_{j}^{s}}\geq \xi_{j}-3\mathfrak{m}_{j}-e^{-\frac{1}{20}N_{j}^{s}},
\end{equation}
\begin{equation}\label{m-5.18}
\mathfrak{m}_{j+1} \leq \mathcal{M}_{j}e^{-5N^{2s-1}N_{j}^{s}}\leq e^{-4N^{2s-1}N_{j}^{s}}.
\end{equation}
By the induction hypothesis,~$\mathfrak{m}_{j}\leq\frac{1}{10}\xi_{j}$ and $\frac{1}{2^{j}}\xi_{j_{0}+1}\leq\xi_{j}$.~Then,~we can use Lemma \ref{lem5.3},~(\ref{m-5.17}-\ref{m-5.18}) and \eqref{11},~to get
$$
  \begin{array}{ll}
    \xi_{j+1}&\geq\frac{7}{10}\xi_{j}-e^{-\frac{1}{20}N_{j}^{s}} \\
    &\geq\frac{7}{10\times 2^{j}}\xi_{j_{0}+1}-e^{-\frac{1}{20}N_{j}^{s}}\\
    &=\frac{1}{5\times 2^{j}}\xi_{j_{0}+1}+\frac{1}{ 2^{j+1}}\xi_{j_{0}+1}-e^{-\frac{1}{20}N_{j}^{s}}\\
    &\geq\frac{1}{5\times 2^{j}}e^{-\frac{1}{1000}N^{s-1}N_{j_{0}}^{s}}-e^{-\frac{1}{20}N_{j}^{s}}+\frac{1}{ 2^{j+1}}\xi_{j_{0}+1}\\
    &\geq\frac{1}{10\times 2^{j}}e^{-\frac{1}{1000}N^{s-1}N_{j_{0}}^{s}}+\frac{1}{ 2^{j+1}}\xi_{j_{0}+1}\\
    &\geq10e^{-4N^{2s-1}N_{j_{0}}^{s}}+\frac{1}{ 2^{j+1}}\xi_{j_{0}+1}\\
    &\geq10\mathfrak{m}_{j+1}+\frac{1}{ 2^{j+1}}\xi_{j_{0}+1}.
  \end{array}
$$
\end{pf}

\subsection{Last Proof Of Theorem \ref{th-1.2} }
 In lemma \ref{lem4.1},~by $B_{j}$  the $cocycle~(\alpha,Ae^{F(\theta)})$ can be conjugated to  $cocycle~ (\alpha,A_{j}e^{F_{j}+\Sigma_{p=j+1}^{\infty}\lambda v_{p}W_{j}})$,~where $W_{j}=Ad(\widetilde{B}_{j}).W $ and the rotation number $\tilde{\rho}_{j}$ of
$cocycle~(\alpha,A_{j}e^{F_{j}+\Sigma_{p=j+1}^{\infty}\lambda v_{p}W_{j}})$ is
\begin{equation}\label{5.14}
  \tilde{\rho}_{j}=\frac{1}{2}\langle k_{J}-\tilde{k}_{j},\alpha\rangle\quad(mod\quad\mathbb{Z}).
\end{equation}
Lemma \ref{lem4.1} guarantees that $cocycle~(\alpha,A_{j}e^{F_{j}+\Sigma_{j+1}^{\infty}\lambda v_{j}W_{j}})$ converges to a constant $cocycle~(\alpha,A_{\infty})$.~In the following,~we will prove that the rotation number of $A_{j}$ is zero and $\|A_{j}\|\geq\delta>0$ for all $j>J+1$ which implies that $A_{\infty}$ is hyperbolic or parabolic.

\begin{Lemma}\label{lem5.6}
 \begin{equation}\label{5.15}
  \|A_{J+1}-I\|\geq |\lambda|e^{-\frac{6}{5}N_{J+1}^{s}},\quad \tilde{\rho}_{J+1}=0.
 \end{equation}
\end{Lemma}
\begin{pf}
Recall $\rho_{j}=rot(\alpha,A_{j})\quad(mod\quad\mathbb{Z})$,~By Lemma \ref{lem2.1} and Lemma \ref{lem4.1},
\begin{equation}\label{5.16}
  \begin{array}{ll}
  |\rho-\tilde{\rho}_{j}|&\leq15(\sup_{\theta\in\mathbb{T}^{d}}\|F_{j}+\Sigma_{p=j}^{\infty}v_{p}W_{p}\|)^{\frac{1}{2}}\\
  &\leq 15(\sigma e^{-N_{j+1}^{s}}+\sigma e^{\frac{1}{20}N_{j}^{s}\Sigma_{p=j}^{\infty}e^{-\frac{9}{10}|k_{p}|^{s}}})^{\frac{1}{2}} \\
  &\leq e^{-\frac{1}{4}N_{j}^{s}},
   \end{array}
\end{equation}
which implies
\begin{equation}\label{5.17}
  |2\rho_{j}-\langle k_{J}-\tilde{k}_{j},\omega\rangle|\leq2e^{-\frac{1}{4}N_{j}^{s}}.
\end{equation}

In case that $J=1$,~we have $\tilde{\rho}_{1}=\frac{1}{2}\langle k_{1},\alpha\rangle\quad(mod\quad\mathbb{Z})$.~By \eqref{5.17},
$$A_{1}=A\in\mathcal{RS}(40N^{s}N_{2},e^{-\frac{1}{50}N_{2}^{s}}),\quad \breve{k}_{1}=k.$$

In case that $J>1$,~we first prove that either $A_{J-1}\in\mathcal{RS}(40N^{s}N_{J},e^{-\frac{1}{50}N_{J-1}^{s}})$ or $A_{J}\in\mathcal{RS}(40N^{s}N_{J},e^{-\frac{1}{50}N_{J}^{s}})$.~In fact,~if $|k_{J}-\tilde{k}_{J-1}|\leq40N^{s}N_{J}$,~then by \eqref{5.17}
$$A_{J-1}\in\mathcal{RS}(40N^{s}N_{J},e^{-\frac{1}{50}N_{J-1}^{s}}),~\quad \breve{k}_{J-1}=k_{J}-\tilde{k}_{J-1}.$$
Otherwise,~if $|k_{J}-\tilde{k}_{J-1}|>40N^{s}N_{J}$,~then by \eqref{5.16} and \eqref{5.17}
$$\begin{array}{ll}
  |2\rho_{J-1}-\langle k,\alpha\rangle|&\geq |\langle k-(k_{J}-\tilde{k}_{J-1}),\alpha\rangle|-|2\rho_{J-1}-\langle k-(k_{J}-\tilde{k}_{J-1}),\alpha\rangle|\\
  &\geq2e^{-\frac{1}{50}N_{J}^{s}}-2e^{-\frac{1}{4}N_{J}^{s}} \\
  &>e^{-\frac{1}{50}N_{J}^{s}},\quad k\in\mathbb{Z}^{d}\backslash\{0\},\quad |k|\leq40N^{s}N_{J},
\end{array}$$
which implies that $A_{J-1}\in\mathcal{NR}(40N^{s}N_{J},e^{-\frac{1}{50}N_{J-1}^{s}})$.~Thus $\tilde{k}_{J}=\tilde{k}_{J-1}$ and by
\ref{lem4.1}
$$|k_{J}-\tilde{k}_{J}|=|k_{J}-\tilde{k}_{J-1}|\leq N_{J+1}+41N^{s}N_{J}\leq40N^{s}N_{J+1}.$$
Then,~we have
$$A_{J}\in\mathcal{RS}(40N^{s}N_{J+1},e^{-\frac{1}{50}N_{J}^{s}}),\quad \breve{k}_{J}=k_{J}-\tilde{k}_{J}.$$

Now,~we divide the proof lemma \ref{lem5.6} into two different cases:\\
(1)$A_{J-1}\in\mathcal{RS}(40N^{s}N_{J},e^{-\frac{1}{50}N_{J-1}^{s}})$ with $J>1$,~Then by \eqref{3.6} $A_{J-1}=\tilde{A}_{J-1}^{'}e^{H
_{J-1}}=e^{A_{J-1}^{''}}.$~
Using Lemma \ref{lem2.1},~\eqref{5.14} and Lemma \ref{lem4.1},~we have
$$k_{J}=\tilde{k}_{J},\quad \tilde{\rho}_{J}=0,\quad |\breve{k}_{J-1}|\geq\frac{999}{1000}|k_{J}|,\quad
A_{J}\in\mathcal{NR}(40N^{s}N_{J+1},e^{-\frac{1}{50}N_{J}^{s}}),$$
Then $A_{J}=A_{J-1}e^{\langle \lambda v_{J}W_{J}\rangle}$.~
Note that $\mathcal{K}\cap\mathcal{Z}_{J}\neq\emptyset\Rightarrow\mathcal{K}\cap\mathcal{Z}_{J-1}=\emptyset$.~
Let us write $$A_{J-1}=\left(
                               \begin{array}{cc}
                                 ia_{J} & b_{J} \\
                                 \overline{b_{J}} & -ia_{J} \\
                               \end{array}
                             \right),$$
where $W_{J}=Ad(\tilde{B}_{J}).W_{E}$ which is written as (note that in this case $\tilde{k}_{J}=k_{J}$)
$$W_{J}=\left(
          \begin{array}{cc}
            iu_{J} & w_{J}e^{-i\langle \tilde{k}_{J},\cdot\rangle} \\
            \overline{w_{J}}e^{i\langle \tilde{k}_{J},\cdot\rangle} & -iu_{J} \\
          \end{array}
        \right)=\left(
          \begin{array}{cc}
            iu_{J} & w_{J}e^{-i\langle k_{J},\cdot\rangle} \\
            \overline{w_{J}}e^{i\langle k_{J},\cdot\rangle} & -iu_{J} \\
          \end{array}
        \right)
$$
By \eqref{4.37} in lemma \ref{lem4.1} and \eqref{11}(note that $s\in(0,\frac{1}{2})$),

\begin{equation}\label{est-b}
\begin{array}{ll}
  |b_{J}|&\leq |\lambda|e^{-\frac{3}{4}|\breve{k}_{J-1}|h_{J-2}}\leq|\lambda| e^{-\frac{3}{4}\times\frac{999}{1000}|k_{J}|h_{J-2}} \\
  &\leq |\lambda|e^{-\frac{1}{20}|k_{J}|^{s}(\frac{|k_{J}|}{N_{J-1}})^{1-s}}\leq |\lambda|e^{-\frac{1}{20}|k_{J}|^{s}N^{1-s}} \\
  &\leq |\lambda|e^{-|k_{J}|^{s}}.
\end{array}
\end{equation}

In view of $v_{J}W_{J}=e^{-|k_{J}|^{s}}(e^{i\langle k_{J},\theta\rangle}+e^{-i\langle k_{J},\theta\rangle})W_{J}$,
by Lemma \ref{lem5.5} (note that $J>J-1>j_{0}$ and $k_{J}=\tilde{k}_{J}$),
$$\begin{array}{ll}
  \|A_{J}e^{\langle\lambda  v_{J}W_{J}\rangle}-I\|&\geq |\lambda|e^{-|k_{J}|^{s}}(|\langle w_{J}\rangle|-|\hat{w}_{J}(2k_{J})|)-|b_{J}|\\
  &\geq |\lambda|e^{-|k_{J}|^{s}}(\xi_{J}-m_{J})-|\lambda|e^{-|k_{J}|^{s}}\\
  &\geq|\lambda|\frac{9}{10} e^{-|k_{J}|^{s}}\xi_{J}-|\lambda| e^{-|k_{J}|^{s}}\\
  &\geq|\lambda|\frac{1}{2^{J+2}}e^{-\frac{1}{1000}N^{s-1}N_{j_{0}}^{s}}e^{-|k_{J}|^{s}}- |\lambda| e^{-|k_{J}|^{s}}\\
  &\geq |\lambda|\frac{1}{2^{J+3}}e^{-\frac{1}{1000}N^{s-1}N_{j_{0}}^{s}}e^{-|k_{J}|^{s}}\\
  &\geq |\lambda|e^{-\frac{11}{10}|k_{J}|^{s}}
\end{array}$$
Now,~by \eqref{4.14} in Lemma \ref{lem4.1}
$$\|A_{J+1}-I\|\geq |\lambda|e^{-\frac{11}{10}|k_{J}|^{s}}-|\lambda| e^{-\frac{8}{5}|k_{J}|^{s}}\geq |\lambda| e^{-\frac{7}{5}|k_{J}|^{s}}\geq |\lambda|e^{-\frac{6}{5}N_{J+1}^{s}}.$$
By \ref{lem2.2}, we get $\tilde{\rho}_{J+1}=\tilde{\rho}_{J}=0$
(2)$A_{J}\in \mathcal{RS}(40N^{s}N_{J+1},e^{-\frac{1}{50}N_{J}^{s}})$.~There is $P_{J}\in SU(1,1)$,~such that
$$Ad(P^{-1}_{J}).A_{J}=\left(
                        \begin{array}{cc}
                          e^{i\rho_{J}} & 0 \\
                          0 & e^{-i\rho_{J}} \\
                        \end{array}
                      \right)~~~~~(0\neq\rho_{J}\in\mathbb{R}),
$$

$$\|P_{J}\|^{2}\leq\frac{1}{1+|\rho_{J}|}e^{\frac{1}{1500}N^{s-1}N_{J}^{s}}.$$
Define
$$\widetilde{W}_{J}\triangleq Ad(P^{-1}_{J}\tilde{B}_{J}).W=Ad(P^{-1}_{J}).W_{J}\triangleq\left(
                                                                                    \begin{array}{cc}
                                                                  i\tilde{u}_{J} & \tilde{w}_{J}e^{-i\langle\tilde{k}_{J},\cdot\rangle} \\
                                                     \overline{\widetilde{w}_{J}}e^{i\langle\tilde{k}_{J},\cdot\rangle}& -i\tilde{u}_{J} \\
                                                                                    \end{array}
                                                                                  \right),
$$
and then $Ad(P^{-1}_{J}\tilde{B}_{J}).(v_{J}W)=v_{J}\widetilde{W}_{J}.$~Write
$$v_{J}\widetilde{W}_{J}=\left(
                           \begin{array}{cc}
                             if_{J}^{+} & g_{J}^{+} \\
                             \overline{g_{J}^{+}} & -if_{J}^{+} \\
                           \end{array}
                         \right)
$$
In view of $v_{J}\widetilde{W}_{J}=e^{-|k_{J}|^{s}}(e^{i\langle k_{J},\theta\rangle}+e^{-i\langle k_{J},\theta\rangle})\widetilde{W}_{J},$~
we have
$$\hat{g}_{J}^{+}(k_{J}-\tilde{k}_{J})=e^{-|k_{J}|^{s}}(\langle \widetilde{W}_{J}\rangle+\widehat{\widetilde{W}}_{J}(2k_{J})),$$
and then
\begin{equation}\label{5.18}
  |\hat{g}_{J}^{+}(k_{J}-\tilde{k}_{J})|\geq e^{-|k_{J}|^{s}}(|\langle \widetilde{W}_{J}\rangle|-|\widehat{\widetilde{W}}_{J}(2k_{J})|).
\end{equation}
When $J>j_{0}\geq1,$~by \eqref{5.18},~Lemma \ref{lem5.5}
$$\begin{array}{ll}
  |\hat{g}_{J}^{+}(k_{J}-\tilde{k}_{J})|&\geq e^{-|k_{J}|^{s}}(|\langle \widetilde{W}_{J}\rangle|-|\widehat{\widetilde{W}}_{J}(2k_{J})|) \\
  &\geq e^{-|k_{J}|^{s}}\{\frac{\|P\|^{2}+1}{2}(\xi_{J}-3m_{J}-)\|P\|^{2}m_{J}\}\\
  &\geq\frac{1}{4} e^{-|k_{J}|^{s}}\xi_{J}\geq\frac{1}{2^{J+3}}e^{-|k_{J}|^{s}}e^{-\frac{1}{1000}N^{s-1}N_{j_{0}}^{s}}\\
  &\geq e^{-\frac{11}{10}|k_{J}|^{s}}
\end{array}$$
As for $J=j_{0}$,~by \eqref{5.18},~Lemma\ref{lem5.2}
$$|\langle\tilde{w}_{j_{0}}\rangle|\geq\frac{1}{4|\rho_{j_{0}}|\|P_{j_{0}}\|^{2}}\geq e^{-\frac{1}{1000}N^{s-1}N_{j_{0}}^{s}},$$
$$|\tilde{w}_{j_{0}}|_{h_{j_{0}}}\leq \|P_{j_{0}}\|^{2}|w_{j_{0}}|_{h_{j_{0}}}\leq e^{\frac{1}{1000}N^{s-1}N_{j_{0}}^{s}},$$
and then (note that $|k_{j_{0}}|h_{j_{0}}\geq\frac{1}{10}N_{j_{0}}N_{j_{0}+1}^{s-1}=\frac{1}{10}N_{j_{0}}N_{j_{0}}^{s}$)
$$\begin{array}{ll}
  |\hat{g}_{j_{0}}^{+}(k_{j_{0}}-\tilde{k}_{j_{0}})|&\geq e^{-|k_{j_{0}}|^{s}}(|\langle \widetilde{W}_{j_{0}}\rangle|-| \widetilde{W}_{j_{0}}|_{h_{j_{0}}}e^{-2|k_{j_{0}}|h_{j_{0}}}) \\
  &\geq e^{-|k_{j_{0}}|^{s}}(e^{-\frac{1}{1000}N^{s-1}N_{j_{0}}^{s}}-e^{\frac{1}{1000}N^{s-1}N_{j_{0}}^{s}}e^{-\frac{1}{5}N^{s-1}N_{j_{0}}^{s}})\\
  &\geq\frac{1}{2} e^{-|k_{j_{0}}|^{s}}e^{-\frac{1}{1000}N^{s-1}N_{j_{0}}^{s}}\\
  &\geq e^{-\frac{11}{10}|k_{j_{0}}|^{s}}.
\end{array}$$
In any case,~by \eqref{3.6},we get
$A_{J}=\tilde{A}_{J}^{'}e^{\langle \lambda v_{J}\widetilde{W}_{J}\rangle}$,~then
$$\|A_{J}-I\|\geq|\lambda||\hat{g}_{J}^{+}(k_{J}-\tilde{k}_{J})|\geq|\lambda| e^{-\frac{11}{10}|k_{J}|^{s}}.$$
By \eqref{4.14} in Lemma \ref{lem4.1},~we get
$$\begin{array}{ll}
  \|A_{J+1}-I\|&\geq|\lambda| e^{-\frac{11}{10}|k_{J}|^{s}}-|\lambda| e^{-\frac{7}{5}|k_{J}|^{s}} \\
  &\geq |\lambda|e^{-\frac{6}{5}|k_{J}|^{s}} \\
  &\geq|\lambda| e^{-\frac{6}{5}N_{J+1}^{s}}.
\end{array}$$
By $\breve{k}=k_{J}-\tilde{k}_{J}$, \eqref{5.14} and Lemma \ref{lem2.2}, we get $\tilde{\rho}_{J+1}=0$.
\end{pf}

\begin{Lemma}\label{lem5.7}
$A_{j}\in \mathcal{NR}(40N^{s}N_{j+1},e^{-\frac{1}{50}N_{j}^{s}})$ and
\begin{equation}\label{5.19}
  |\lambda|e^{-2N_{J+1}^{s}}\leq\|A_{j}-I\|\leq|\lambda| e^{-\frac{1}{10}N_{J+1}^{s}},\quad \tilde{\rho}_{j}=0\quad for ~j\geq J+2
\end{equation}
\end{Lemma}

\begin{pf}
We prove inductively that for all $j\geq J+2$
$$A_{j}\in \mathcal{NR}(40N^{s}N_{j+1},e^{-\frac{1}{50}N_{j}^{s}}),\quad \|A_{j}-I\|\geq(1+\frac{1}{2^{j}}e^{-2N_{J+1}^{s}}),\quad \tilde{\rho}_{j}=0$$
When $j=J+2$,~$\mathcal{K}\cap\mathcal{Z}_{J}\neq\emptyset$~implies that $\mathcal{K}\cap\mathcal{Z}_{J}=\emptyset$.~By \eqref{4.14} in Lemma \ref{lem4.1},~Lemma \ref{lem5.6} and \eqref{11},~we get
$$\|A_{J+2}-I\|\geq\|A_{J+1}-I\|-|\lambda| e^{-\frac{7}{5}|k_{J}|^{s}}\geq |\lambda|e^{-\frac{7}{5}|k_{J}|^{s}}\geq |\lambda|(1+\frac{1}{2^{J+2}}e^{-2N_{J+1}^{s}}).$$
By Lemma \ref{lem2.2} and \eqref{5.16},~we have $\tilde{\rho}_{J+2}=0$ and $|\rho_{J+2}|<e^{-\frac{1}{4}N_{J+2}^{s}}.$~It follows that
$A_{j+2}\in \mathcal{NR}(40N^{s}N_{J+3},e^{-\frac{1}{50}N_{J+2}^{s}})$.
We inductively assume that the desired conclusion holds for $j\geq J+2$,~we now verify it for $j+1$.~In fact,~by \eqref{4.14} and \eqref{11},
$$
\begin{array}{ll}
 \|A_{j+1}-I\|&\geq\|A_{j}-I\|-|\lambda| e^{-\frac{3}{5}N_{j}^{s}} \\
  &\geq |\lambda|(1+\frac{1}{2^{j}}e^{-2N_{J+1}^{s}})- |\lambda| e^{-\frac{3}{5}N_{j}^{s}} \\
  &\geq |\lambda|(1+\frac{1}{2^{j+1}}e^{-2N_{J+1}^{s}}).
\end{array}
$$
By Lemma \ref{lem2.2},~$\tilde{\rho}_{j+1}=\tilde{\rho}_{j}=0$.~Then,~for $|\rho_{j+1}|<e^{-\frac{1}{4}N_{J+1}^{s}}$ and \eqref{5.16},~we bave
$A_{j+1}\in \mathcal{NR}(40N^{s}N_{j+2},e^{-\frac{1}{50}N_{j+1}^{s}})$.

For the upper bounds of $\|A_{j}-I\|$, similar to Lemma \ref{lem5.6}, we are going to prove it in two cases:

(1)$A_{J-1}\in \mathcal{RS}(40N^{s}N_{j+1},e^{-\frac{1}{50}N_{j}^{s}})$. By \eqref{est-b} and \eqref{4.15}, we have
\begin{equation}\label{est-bj}
  \begin{array}{ll}
  |b_{j}|&\leq|b_{J}|+\sum_{p=J}^{j-1}|b_{p+1}-b_{p}|\\
  &\leq |\lambda|e^{-|k_{J}|^{s}}+|\lambda| e^{-\frac{2}{5}N_{J+1}^{s}}\\
  &\leq |\lambda|e^{-\frac{1}{5}N_{J+1}^{s}}
  \end{array}
\end{equation}

(2)$A_{J-1}\in \mathcal{NR}(40N^{s}N_{j+1},e^{-\frac{1}{50}N_{j}^{s}})$. By \eqref{4.1} and \eqref{4.33}, we have
$$|P_{J-1}Ad(\widetilde{B}_{J-1}).(\lambda v_{J-1}W)P^{-1}_{J-1}|_{h_{J-1}}\leq 2|\lambda|e^{-\frac{1}{5}N_{J+1}^{s}},$$
then the inequality in right hand of  \eqref{5.7} is a direct result of \eqref{4.37} in Lemma \ref{lem4.1}.
\end{pf}

Now,~we have proved that $A_{j}\in \mathcal{NR}(40N^{s}N_{j+1},e^{-\frac{1}{50}N_{j}^{s}})$ for all $j\geq J+2$,~thus $B_{j}$ always satisfies $|B_{j}-I|_{h_{j}}\leq |\lambda|e^{-\frac{3}{5}N_{j}^{s}}$ (by \eqref{4.15}). Then, for all $j\geq J+2$, we have
\begin{equation}\label{jb}
  \begin{array}{ll}
    |\widetilde{B}_{j}|_{h_{j-1}}&=|B_{j-1}B_{j-2}B_{J+1}\widetilde{B}_{J+1}|_{h_{j-1}}\\
    &=(\prod_{p=J+1}^{j-1}|B_{p}|_{h_{p}})|\widetilde{B}_{J+1}|_{h_{J+1}}\\
    &\leq\prod_{p=J+1}^{j-1}(1+|\lambda| e^{-\frac{3}{5}N_{p}^{s}})e^{\frac{1}{40}N_{j+1}^{s}}\\
    &\leq 2e^{\frac{1}{40}N_{J+1}^{s}},
  \end{array}
\end{equation}
and
\begin{equation}\label{bb}
\begin{array}{ll}
  \|\widetilde{B}_{j+1}-\widetilde{B}_{j}\|_{C^{0}}&\leq|\widetilde{B}_{j+1}-\widetilde{B}_{j}|_{h_{j}}\\
  &\leq |\widetilde{B}_{j}|_{h_{j}}|B_{j}-I|_{h_{j}}\\
  &\leq 2|\lambda| e^{-\frac{23}{40}N_{J+1}^{s}}.
  \end{array}
\end{equation}
Then $\widetilde{B}_{j}$ is convergent in $C^{\infty}$ topology.

Set $\widetilde{B}_{\infty}=\lim_{j\rightarrow\infty}\widetilde{B}_{j}$, then $\widetilde{B}_{\infty}$ conjugates the $Cocycle~(\alpha,Ae^{F(\theta)})$ with $v$ defined in \eqref{1.3} to $Cocycle~(\alpha,A_{\infty})$ ~$(A_{\infty}=\lim_{j\rightarrow\infty}A_{j})$.
Its fibered rotation number is zero because $\lim_{j\rightarrow\infty}\tilde{\rho}_{j}=0$.~Thus the key point is to prove $A_{\infty}\neq I$.
~This is an obvious consequence of \eqref{5.19},~for
$$\|A_{\infty}-I\|=\lim_{j\rightarrow\infty\|A_{j}-I\|}\geq |\lambda|e^{-2N_{J+1}^{s}}>0.$$
We then arrive at the conclusion that the $Cocycle~(\alpha,A_{\infty})$ is either uniformly hyperbolic or parabolic. Then by Lemma \ref{Lem-A2}, there exists $R_{\phi}$ such that $R_{-\phi}A_{\infty}R_{\phi}=A$ where $A=\left(
                                                                                           \begin{array}{cc}
                                                                                             1 & \zeta \\
                                                                                             0 & 1 \\
                                                                                           \end{array}
                                                                                         \right)$, and from \eqref{5.19} we have
\begin{equation}\label{c}
 |\lambda| e^{-2N_{J+1}^{s}}\leq|\zeta|\leq |\lambda|e^{-\frac{1}{6}N_{J+1}^{s}}.
\end{equation}

\hspace*{\fill}\
\section{ PROOF OF THEOREM \ref{main-th} }
Next, under the assumption of Theorem \ref{th-1.2}, we will proof the Theorem \ref{main-th}, The theorem to be proved in this section is given below.

\begin{Theorem}\label{th-6.1}
Under the assumption of Theorem \ref{th-1.2},~for any $k\in\mathcal{K}$,~the spectral gap $I_{k}(\lambda v)$ has the following size:
$$|\lambda|^{2}e^{-\frac{13}{6}|k|^{2s}}\leq|I_{k}(\lambda v)|\leq \sqrt{|\lambda|}e^{-\frac{3}{20}|k|^{s}}.$$
\end{Theorem}

Before we start to prove it we will focus on gap estimates for  quasi-periodic Schr\"odinger operator on $l^{2}(\mathbb{Z})$:
$$(H_{v,\alpha,\theta}u)_{n}=u_{n+1}+u_{n-1}+\lambda v(\theta+n\alpha)u_{n},$$
with $\alpha\in\mathbb{T}^{d}$ such that $(1,\alpha)$ is rationally independent,~and  Gevery potential $v=\sum_{k\in\mathcal{K}}e^{-|k|^{s}}cos\langle k,\theta\rangle$ .~We will estimate the size of the spectral gap $I_{k}(\lambda v)=(E_{k}^{-},E_{k}^{+})$ via Moser-P\"oschel argument \cite{MP84} at its edge points.

 Then from Lemma \ref{lem4.1} , we find that the cocycle $(\alpha,S_{E_{k}}^{\lambda v})$ is reduced to
 $(\alpha, Ae^{F_{j+1}(\theta)+Ad(\widetilde{B}_{j}).(\sum_{p=j+1}^{\infty}(\lambda v_{p}W_{E}))})$ at the step $j$ with $\widetilde{B}_{j}\in C_{h_{j}}^{\omega}(\mathbb{T}^{d},PSL(2,\mathbb{R}))$ for some $0<h_{j}<1$ ~such that
$$\widetilde{B}_{j}(\cdot+\alpha)^{-1}S_{E_{k_{j}}}^{\lambda v}\widetilde{B}_{j}(\cdot)=Ae^{F_{j+1}(\theta)+Ad(\widetilde{B}_{j}).(\sum_{p=j+1}^{\infty}(\lambda v_{p}W_{E }))},$$
with
$ A=\left(
     \begin{array}{cc}
       1 & \zeta \\
       0 & 1 \\
     \end{array}
   \right)$
~  and
 \begin{equation}\label{6.1}
 |\widetilde{B}_{j}|_{h_{j}}\leq e^{\frac{1}{40}N_{j}^{s}}, \quad|\lambda| e^{-2N_{J+1}^{s}}\leq|\zeta|\leq |\lambda|e^{-\frac{1}{6}N_{J+1}^{s}}.
  \end{equation}
 Moreover, we can deduce from \eqref{bb} that
 \begin{equation}\label{inbb}
   \|\widetilde{B}_{\infty}-\widetilde{B}_{j}\|_{C^{0}}\leq \sum_{p=j}^{\infty}\|\widetilde{B}_{j+1}-\widetilde{B}_{j}\|_{C^{0}}\leq 2|\lambda| e^{-\frac{1}{2}N_{J+1}^{s}}.
 \end{equation}
 And let $Ae^{F_{j+1}(\theta)+Ad(\widetilde{B}_{j}).(\sum_{p=j+1}^{\infty}(\lambda v_{p}W_{E}))}:=A+M_j$, by \eqref{4.1} and $\|Ad(\widetilde{B}_{j}).(\lambda v_{p}W_{E})\|_{C^{0}}\leq |\lambda|e^{-\frac{17}{20}|k|_{j}^{s}}$ we get
 \begin{equation}\label{est-M}
 \quad \|M_j\|_{C^{0}}\leq\frac{3}{2}\|A\|(|F_{j+1}|_{h_{j}}+\|Ad(\widetilde{B}_{p}).(\lambda v_{p}W_{E})\|_{C^{0}})\leq 3e^{-\frac{4}{5}|k|_{j}^{s}}.
 \end{equation}

\subsection{Moser-P\"oschel Argument}
For any $0<\delta<1$,~we can calculate
$$\widetilde{B}_{j}(\cdot+\alpha)^{-1}S_{E_{k_{j}-\delta}^{+}}^{\lambda v}\widetilde{B}_{j}(\cdot)=A+M_j-\delta P(\cdot)$$
where
$$\widetilde{B}_{j}=\left(
                      \begin{array}{cc}
                        B_{j,11} & B_{j,12}\\
                        B_{j,21} & B_{j,22} \\
                      \end{array}
                    \right)
$$
$$P(\cdot):=\left(
              \begin{array}{cc}
                B_{j,11}(\cdot)B_{j,12}(\cdot)-\zeta B_{j,11}^{2}(\cdot) & -\zeta B_{j,11}(\cdot)B_{j,12}(\cdot)+B_{j,12}^{2}(\cdot) \\
                -B_{j,11}^{2}(\cdot) & -B_{j,11}(\cdot)B_{j,12}(\cdot) \\
              \end{array}
            \right).
$$
And,
\begin{equation}\label{6.0}
  |P|_{h^{'}_{j}}\leq(1+\zeta )|\widetilde{B}_{j}|^{2}_{h^{'}_{j}},\quad h^{'}_{j}\leq \frac{1}{90}h_j .
\end{equation}

 Then, We find that when the energy $E$ move from the right end of the gap $E_{k}^{+}$ to $E_{k}^{+}-\delta$,~ the spectral gap of the other edge point is  according to the variation of the rotation number $\rho(\alpha,\widetilde{B}_{j}(\cdot+\alpha)^{-1}S_{E_{k}^{+}-\delta}^{\lambda v}\widetilde{B}_j(\cdot))$.~And the rotation number of the constant cocycle $(\alpha,A)$ is zero since $A$ is parabolic.

In the following,~we first apply one standard KAM step to the cocycle $(\alpha,A+M_j-\delta P(\cdot))$,~which serves as the starting point of our estimation of the size of the gap. We denote $[\cdot]$ the average of a quasi-periodic function.

\begin{Lemma}\cite[Lemma6.1]{LYZZ}\label{lem-6.1}
Given $\alpha\in DC(\gamma,\tau)$, $\tau>d-1$ let
  $c_{\tau}:=2^{8}\Gamma(4\tau+2)$ , if~~~$0<\delta<c_{\tau}^{-1}\gamma^{3}h_{j}^{'4\tau+1}|\widetilde{B}_j|_{h^{'}_{j}}^{-2}$,~then there exist $\tilde{B}\in C^{\omega}_{\tilde{h}_j}(\mathbb{T}^{d},SL(2,\mathbb{R}))$ and $P_{1}\in C^{\omega}_{\tilde{h}_j}(\mathbb{T}^{d},gl(2,\mathbb{R}))$ ($\tilde{h}_j<\frac{1}{2}h^{'}_{j}$)such that
\begin{equation}\label{6.2}
  \tilde{B}(\cdot+\alpha)^{-1}(A+M_j-\delta P(\cdot))\tilde{B}=e^{b_{0}-\delta b_{1}}+\delta^{2}P_{1}(\cdot),
\end{equation}
where $b_{0}:=\left(
                \begin{array}{cc}
                  0 & \zeta  \\
                  0 & 0 \\
                \end{array}
              \right)$ and $$b_{1}:=\left(
                                     \begin{array}{cc}
                                       [B_{11}B_{12}]-\frac{\zeta }{2}[B_{11}^{2}] & -\zeta [B_{11}B_{12}]+[B_{12}^{2}] \\
                                       -[B_{11}^{2}] & -[B_{11}B_{12}]+\frac{\zeta }{2}[B_{11}^{2}] \\
                                     \end{array}
                                   \right)$$
with the estimates
\begin{equation}\label{6.3}
\begin{array}{ll}
 & |\widetilde{B}-Id|_{\tilde{h}_j}\leq2|Y|_{\tilde{h}_j}\leq
 2c_{\tau}\gamma^{-3}\delta|\widetilde{B}_j|_{\tilde{h}}^{2}\tilde{h}_j^{-(4\tau+1)},\\
 & \|P_{1}\|_{C^{0}}\leq61c_{\tau}^{2}\gamma^{-6}|\widetilde{B}_j|_{h'_j}^{4}\tilde{h}_j^{-2(4\tau+1)}+\delta^{-1}\zeta^{2}|\widetilde{B}_j|_{h'_j}^{2}.
  \end{array}
\end{equation}
\end{Lemma}

\begin{pf}
 Let $G:=-\delta A^{-1}P$, then $G\in C^{\omega}(\mathbb{T}^{d},sl(2,\mathbb{R}))$. By a standard KAM step, we can solve the linearized cohomological equation
\begin{equation*}
  -Y(\theta+\alpha)(A+M_j)+(A+M_j)Y(\theta)=(A+M_j)(G(\theta)-[G]).
\end{equation*}
Compare the Fourier coefficients of two sides, and by the polynomial decay of Fourier coefficient $\widehat{G}(k)$, we have
\begin{equation*}
\begin{array}{ll}
|Y|_{\tilde{h}_j}&\leq 10\delta |P|_{\tilde{h}_j}\sum_{k\in\mathbb{Z}^{d}}\frac{e^{-h'_j|k|}}{|e^{i\langle k,\alpha\rangle}-1|^{3}}e^{\tilde{h_j}|k|}\\
&\leq 20\gamma^{-3}\delta|P|_{h'_j}\sum_{k\in\mathbb{Z}^{d}}|k|^{3\tau}e^{-(h'_j-\tilde{h}_j)|k|}\\
&\leq 40\gamma^{-3}\delta|P|_{h'_j}\int_{0}^{+\infty}x^{d-1}x^{3\tau}e^{-(h'_j-\tilde{h}_j)x}dx.
\end{array}
\end{equation*}
and $$\int_{0}^{+\infty}x^{d-1}x^{3\tau}e^{-(h'_j-\tilde{h}_j)x}dx\leq2+\int_{0}^{+\infty}x^{4\tau}e^{-(h'_j-\tilde{h}_j)x}dx\leq
2\Gamma(4\tau+2)\cdot\tilde{h}_j^{-(4\tau+1)},$$
so
\begin{equation}\label{est-Y}
|Y|_{\tilde{h}_j}\leq\frac{1}{2}c_{\tau}\gamma^{-3}\delta|P|_{h'_j}\tilde{h}_j^{-(4\tau+1)}
\end{equation}
Let $\widetilde{B}:=e^{Y}$, we have
\begin{equation*}
\widetilde{B}(\theta+\alpha)^{-1}(A+M_j-\delta P(\theta))\widetilde{B}(\theta)=Ae^{[G]}+\widetilde{P}(\theta),
\end{equation*}
with estimate $|\widetilde{B}-Id|_{\tilde{h}_j}\leq2|Y|_{\tilde{h}_j}\leq2c_{\tau}\gamma^{-3}\delta|\widetilde{B}_j|_{\tilde{h}}^{2}\tilde{h}_j^{-(4\tau+1)}$,
where
\begin{equation*}
\begin{array}{ll}
\widetilde{P}(\theta)=&\Sigma_{m+k\geq2}\frac{1}{m!}(-Y(\theta+\alpha))^{m}A\frac{1}{k!}Y(\theta)^{k}-A\Sigma_{k\geq2}\frac{1}{k!}[G]^{k}\\
&-\delta\Sigma_{m+k\geq1}\frac{1}{m!}(-Y(\theta+\alpha))^{m}P(\theta)\frac{1}{k!}Y(\theta)^{k}+M_j\\
&+\Sigma_{m+k\geq1}\frac{1}{m!}(-Y(\theta+\alpha))^{m}M_j(\theta)\frac{1}{k!}Y(\theta)^{k}.
\end{array}
\end{equation*}
Note that $\Sigma_{m+k=l}\frac{l!}{m!k!}=2^{l}$ and  $|G|_{\tilde{h}_j}\leq\delta|P|_{\tilde{h}_j}$, by (\ref{6.1},\ref{est-M},\ref{est-Y}) we get
\begin{equation*}
\begin{array}{ll}
&|\Sigma_{m+k\geq2}\frac{1}{m!}(-Y(\theta+\alpha))^{m}A\frac{1}{k!}Y(\theta)^{k}|_{\tilde{h}_j}
\leq2c_{\tau}^{2}\gamma^{-6}\delta^{2}|P|_{h'_j}^{2}\tilde{h}_j^{-2(4\tau+1)},\\
&|\delta\Sigma_{m+k\geq1}\frac{1}{m!}(-Y(\theta+\alpha))^{m}P(\theta)\frac{1}{k!}Y(\theta)^{k}|_{\tilde{h}_j}
\leq2c_{\tau}\gamma^{-3}\delta^{2}|P|_{h'_j}^{2}\tilde{h}_j^{-(4\tau+1)}\\
&|A\Sigma_{k\geq2}\frac{1}{k!}[G]^{k}|_{\tilde{h}_j}\leq4\delta^{2}|P|_{h'_j}^{2}\\
&\|\Sigma_{m+k\geq1}\frac{1}{m!}(-Y(\theta+\alpha))^{m}M_j(\theta)\frac{1}{k!}Y(\theta)^{k}\|_{C^{0}}
\leq6c_{\tau}\gamma^{-3}\delta^{2}e^{-\frac{9}{10}|k|_{j}^{s}}\tilde{h}_j^{-(4\tau+1)}.
\end{array}
\end{equation*}
Hence,
\begin{equation*}
  \|\widetilde{P}(\theta)\|_{C^{0}}\leq8c_{\tau}^{2}\gamma^{-6}\delta^{2}|P|_{h'_j}^{2}\tilde{h}_j^{-2(4\tau+1)}.
\end{equation*}
Then, we define $\widetilde{P}_{1}:=\delta^{-2}\widetilde{P}+\sum_{j\geq2}(j!)^{-1}(-\delta)^{j-2}A[A^{-1}P]^{j} $ such that
\begin{equation*}
  Ae^{[G]}+\widetilde{P}(\theta)=A-\delta[P]+\delta^{2}\widetilde{P}_{1}(\theta),
\end{equation*}
and
\begin{equation}\label{est-p1}
  \|\widetilde{P}_{1}\|_{C^{0}}\leq 8c_{\tau}^{2}\gamma^{-6}|P(\theta)|_{h'_j}^{2}\tilde{h}_j^{-2(4\tau+1)}+2\times\frac{1}{2!}\|A\|^{3}|P(\theta)|_{h'_j}^{2}
\end{equation}
By direct calculation, we can
\begin{equation*}
  A-\delta[P]=Id+(b_{0}-\delta b_{1})-\frac{\delta}{2}(b_{0}b_{1}+b_{1}b_{0}).
\end{equation*}
Since $b_{0}$ is nilpotent, one can check that
\begin{equation*}
  \widetilde{B}(\theta+\alpha)^{-1}(A+M-\delta P(\theta))\widetilde{B}(\theta)=e^{b_{0}-\delta b_{1}}+\delta^{2}P_{1}(\theta),
\end{equation*}
where $P_{1}(\theta)=\widetilde{P}_{1}-\frac{1}{2}b_{1}^{2}-\delta^{-2}\Sigma_{j\geq3}(j!)^{-1}(b_{0}-\delta b_{1})^{j}$ with estimate
\begin{equation}\label{p11}
\begin{array}{ll}
\|P_{1}(\theta)\|_{C^{0}}&\leq\|\widetilde{P}_{1}(\theta)\|_{C^{0}}+\frac{1}{2}|b_{1}|^{2}+2\delta^{-2}\times\frac{1}{3!}|(b_{0}-\delta b_{1})^{3}|\\
&\leq 8c_{\tau}^{2}\gamma^{-6}|P(\theta)|_{h'_j}^{2}\tilde{h}_j^{-2(4\tau+1)}+2\times\frac{1}{2!}\|A\|^{3}|P(\theta)|_{h'_j}^{2}+\frac{1}{2}|P(\theta)|_{h'_j}^{2}\\
&+\frac{2}{3!}\delta^{-2}(\delta^{3}|P|_{h'_j}^{3}+3\zeta\delta^{2}|P(\theta)|_{h'_j}^{2}+\delta\zeta^{2}|P(\theta)|_{h'_j})\\
&\leq16c_{\tau}^{2}\gamma^{-6}|\widetilde{B}_j|_{\tilde{h}}^{4}\tilde{h}_j^{-2(4\tau+1)}+32|\widetilde{B}_j|_{\tilde{h}}^{4}+2|\widetilde{B}_j|_{\tilde{h}}^{4}\\
&+(9|\widetilde{B}_j|_{h'_j}^{4}+2|\widetilde{B}_j|_{h'_j}^{4}+\delta^{-1}\zeta^{2}|\widetilde{B}_j|_{h'_j}^{2})\\
&\leq61c_{\tau}^{2}\gamma^{-6}|\widetilde{B}_j|_{h'_j}^{4}\tilde{h}_j^{-2(4\tau+1)}+\delta^{-1}\zeta^{2}|\widetilde{B}_j|_{h'_j}^{2}
\end{array}
\end{equation}
by (\ref{6.1},\ref{6.0},\ref{est-Y},\ref{est-p1}). ~
 Hence we finish the proof.
\end{pf}

Since $\widetilde{B}$ is homotopic to identity by construction,~we have
$$\rho(\alpha,A+M_{j}-\delta P(\cdot))=\rho(\alpha,e^{b_{0}-\delta b_{1}}+\delta^{2}P_{1}(\cdot)).$$
Let $d(\delta):=det(b_{0}-\delta b_{1})$.~By a direct calculation,~we get
\begin{equation}\label{6.4}
d(\delta)=-\delta[B_{11}^{2}]c+\delta^{2}([B_{11}^{2}][B_{12}^{2}]-[B_{11}B_{12}]^{2}).
\end{equation}

\subsection{The  bounds of spectral gaps.}

With the assistance of Moser-P\"oschel argument and the reducibility of the Schr\"odinger cocycle, we will use the next Lemma \ref{lem-6.2} to prove Theorem \ref{th-6.1}

Then  we recall the following fundamental lemma which was established  in Lemma 6.2 and Lemma 6.3 of \cite{LYZZ}.
~Firstly make some technical preparations:
\begin{Lemma}\cite[Lemma6.2,6.3]{LYZZ}\label{lem-6.2}
For any $B\in C^{\omega}_{h}(\mathbb{T}^{d},PSL(2,\mathbb{R}))$,~$[B_{11}^{2}]\geq(2|B|_{\mathbb{T}^{d}})^{-2}$.

Moreover, For any $\kappa\in(0,\frac{1}{4})$,~and $\zeta\in(0,\frac{1}{2})$ if
\begin{equation}\label{6.7}
  |B|_{h}\zeta ^{\frac{\kappa}{2}}\leq\frac{1}{4},
\end{equation}
and
$$B=\left(
    \begin{array}{cc}
      B_{11} & B_{12} \\
      B_{21} & B_{22} \\
    \end{array}
  \right)
$$
then the following holds:
\begin{equation}\label{6.8}
  0<\frac{[B_{11}^{2}]}{[B_{11}^{2}][B_{12}^{2}]-[B_{11}B_{12}]^{2}}\leq\frac{1}{2}\zeta ^{-\kappa},
\end{equation}
\begin{equation}\label{6.9}
  [B_{11}^{2}][B_{12}^{2}]-[B_{11}B_{12}]^{2}\geq8\zeta ^{2\kappa}.
\end{equation}
\end{Lemma}

Then under the Lemma \ref{lem-6.2}, we proof of Theorem \ref{th-6.1}.
\begin{pf}(Proof of Theorem \ref{th-6.1})
By \eqref{6.1} we have
\begin{equation}\label{6.10}
  |\widetilde{B}_{j}|_{h^{'}_{j}}^{14}\zeta^{\frac{1}{17}}\leq e^{\frac{14}{40}N_{j+1}^{s}}e^{-\frac{3}{26}N_{j+1}^{s}}\leq 10^{-11}c_{\tau}^{-4}\gamma^{12}.
\end{equation}
And, by \eqref{6.4},~the quantity $d(\delta)=det(b_{0}-\delta b_{1})$ satisfies
$$\begin{array}{ll}
    d(\delta)&=-\delta[B_{j,11}^{2}]c+\delta^{2}([B_{j,11}^{2}][B_{j,12}^{2}]-[B_{j,11}B_{j,12}]^{2}) \\
    &=\delta([B_{j,11}^{2}][B_{j,12}^{2}]-[B_{j,11}B_{j,12}]^{2})(\delta-\frac{[B_{j,11}^{2}]\zeta }{[B_{j,11}^{2}][B_{j,12}^{2}]-[B_{j,11}B_{j,12}]^{2}}).
  \end{array}
$$
 Let $\delta_{1}:=\zeta ^{\frac{16}{17}}$, and recall $\zeta$ by \eqref{c}. By \eqref{6.10}, we have
 $$\delta_{1}c_{\tau}\gamma^{-3}|\widetilde{B}_{j}(\theta)|_{h'_{j}}^{2}\leq
 \zeta^{\frac{4}{17}}c_{\tau}\gamma^{-3}|\widetilde{B}_{j}(\theta)|_{h'_{j}}^{\frac{7}{2}}
 \leq10^{-\frac{5}{4}}<1.$$
 For $\alpha\in DC_{d}(\gamma,\tau)$, so $0<\delta_{1}<c_{\tau}^{-1}\gamma^{3}|\widetilde{B}_{j}|_{h'_{j}}^{-2}$.~Hence,~we can apply Lemma \ref{lem-6.1},~the cocycle $(\alpha,A+M-\delta_{1}P(\theta))$ is conjugated to the cocycle $(\alpha,e^{b_{0}-\delta_{1}b_{1}}+\delta_{1}^{2}P_{1})$.~By \eqref{6.10},~one has $|\widetilde{B}_{j}|_{h'_{j}}\zeta ^{\frac{1}{34}}\leq10^{-\frac{5}{2}}c_{\tau}^{-2}\gamma^{6}|\widetilde{B}_{j}|_{h'_{j}}^{-6}\leq\frac{1}{4}$.~Then we can apply Lemma \ref{lem-6.2},~and get
$$\frac{[B_{j,11}^{2}]\zeta }{[B_{j,11}^{2}][B_{j,12}^{2}]-[B_{j,11}B_{j,12}]^{2}}\leq\frac{1}{2}\delta_{1}.$$
Hence,~for $d(\delta_{1})=det(b_{0}-\delta_{1}b_{1})+\frac{1}{4}\delta_{1}^{2}\zeta^{2}[B_{j,11}^{2}]^{2}$,~we have
\begin{equation}\label{6.01}
  d(\delta_{1})\geq \zeta ^{\frac{16}{17}}\cdot8\zeta ^{\frac{2}{17}}\cdot\frac{1}{2}\zeta ^{\frac{16}{17}}=4\zeta ^{2}.
\end{equation}
Following the expressions of $b_{0}$ and $b_{1}$ in Lemma \ref{lem-6.1},~we have
\begin{equation}\label{6.02}
\begin{array}{ll}
  det(b_{0}-\delta_{1}b_{1})&\geq 4\zeta^{2} -\frac{1}{4}\delta_{1}^{2}\zeta^{2}[B_{j,11}^{2}]^{2}\\
  &\geq4\zeta^{2}(1-\frac{1}{16}\zeta^{\frac{32}{17}}|\widetilde{B}_{j}|_{h'_{j}}^{4})\\
  &\geq3\zeta^{2}.
  \end{array}
\end{equation}
In view of Lemma 8.1 in \cite{HY12},~there exists $\mathcal{P}\in SL(2,\mathbb{R})$,~with $\|\mathcal{P}\|\leq2(\frac{\|b_{0}-\delta_{1}b_{1}\|}{\sqrt{d(\delta_{1})}})^{\frac{1}{2}}$ such that
$$\mathcal{P}^{-1}e^{b_{0}-\delta_{1}b_{1}}\mathcal{P}=exp\left(
                                                            \begin{array}{cc}
                                                              0 & \sqrt{det(b_{0}-\delta_{1}b_{1})} \\
                                                              -\sqrt{det(b_{0}-\delta_{1}b_{1})} & 0 \\
                                                            \end{array}
                                                          \right)=:\triangle
.$$
Since $\|b_{0}-\delta_{1}b_{1}\|\leq\zeta+\delta_{1}(1+\zeta)\|\widetilde{B}_{j}\|^{2}_{C^{0}}\leq
\frac{3}{2}\zeta^{\frac{16}{17}}|\widetilde{B}_{j}|^{2}_{h'_{j}}$, combining \eqref{6.01} and \eqref{6.02},~we have
\begin{equation}\label{b0}
\frac{\|b_{0}-\delta_{1}b_{1}\|}{\sqrt{det(b_{0}-\delta_{1}b_{1})}}\leq\frac{2\zeta ^{\frac{16}{17}}|\widetilde{B}_{j}|^{2}_{h'_{j}}}{\sqrt{3}\zeta}\leq|\widetilde{B}_{j}|^{2}_{h'_{j}}\zeta ^{-\frac{1}{17}}.
\end{equation}
Then,~according to Lemma \ref{lem2.1} and Lemma \ref{lem-6.1} with $\mathcal{P}^{-1}(e^{b_{0}-\delta_{1}b_{1}}+\delta_{1}^{2}P_{1})\mathcal{P}=\triangle+\mathcal{P}^{-1}\delta_{1}^{2}P_{1}\mathcal{P}$, by (\ref{p11},\ref{b0}) we have
\begin{equation}\label{x}
\begin{array}{ll}
&|\rho(\alpha,e^{b_{0}-\delta_{1}b_{1}}+\delta_{1}^{2}P_{1})-\sqrt{det(b_{0}-\delta_{1}b_{1})}|\\
&\leq\delta_{1}^{2}\|\mathcal{P}\|^{2}\|P_{1}\|_{C^{0}}\\
&\leq\zeta^{\frac{32}{17}}\times4|\widetilde{B}_{j}|^{2}_{h'_{j}}\zeta^{-\frac{1}{17}}\times
(61c_{\tau}^{2}\gamma^{-6}|\widetilde{B}_{j}|_{h'_{j}}^{4}\tilde{h}_j^{-2(4\tau+1)}+
\zeta^{-\frac{16}{17}}\zeta^{2}|\widetilde{B}_{j}|^{2}_{h'_{j}})\\
&\leq480c_{\tau}^{2}\gamma^{-6}|\widetilde{B}_{j}|_{h'_{j}}^{6}\zeta ^{\frac{31}{17}}\tilde{h}_j^{-2(4\tau+1)}.
\end{array}
\end{equation}
By \eqref{6.10},~we have
$$480c_{\tau}^{2}\gamma^{-6}|\widetilde{B}_{j}|^{6}_{h'_{j}}\zeta ^{\frac{14}{17}}\tilde{h}_j^{-2(4\tau+1)}<1$$
which implies that
\begin{equation*}
\begin{array}{ll}
\rho(\alpha,e^{b_{0}-\delta_{1}b_{1}}+\delta_{1}^{2}P_{1})
&\geq|\rho(\alpha,\triangle)|-|\rho(\alpha,\triangle+\mathcal{P}^{-1}\delta_{1}^{2}P_{1}(\theta)\mathcal{P})-\rho(\alpha,\triangle)|\\
&\geq\sqrt{3}\zeta-480c_{\tau}^{2}\gamma^{-6}|\widetilde{B}_{j}|_{h'_{j}}^{6}\zeta ^{\frac{31}{17}}\tilde{h}_j^{-2(4\tau+1)}\\
&\geq\sqrt{3}\zeta-\zeta>0,
\end{array}
\end{equation*} by (\ref{6.02},\ref{x}).
So by \eqref{c} ,we have $|I_{k}(\lambda v)|\leq\delta_{1}=\zeta^{\frac{16}{17}}\leq  e^{-\frac{3}{20}|k|^{s}}$ with $k$ satisfies (\ref{k-2.7}-\ref{k-2.12}).
This concludes the proof of the upper bound estimates.

Let us now consider the lower bound estimate on the size of the gap.
 Let $\delta_{2}:=\zeta ^{\frac{18}{17}}$.~We are going to show that $|I_{k}(\lambda v)|\geq\delta_{2}$.~We first note that
$$\delta_{2}^{2}|[B_{j,11}^{2}][B_{j,12}^{2}]-[B_{j,11}B_{j,12}]^{2}|\leq2\zeta ^{\frac{36}{17}}|\widetilde{B}_{j}|^{4}_{h'_{j}},$$
and,~by Lemma \ref{lem-6.2},~one has $\delta_{2}[B_{j,11}^{2}]\zeta \geq\frac{1}{4}\zeta ^{\frac{35}{17}}|\widetilde{B}_{j}|^{-2}_{h'_{j}}$.~Thus,~if $\zeta $ is small enough such that
$|\widetilde{B}_{j}|^{6}_{h'_{j}}\zeta ^{\frac{1}{17}}\leq\frac{1}{40}$ (which can be deduced from \eqref{6.10}),~then
\begin{equation*}
\begin{array}{ll}
d(\delta_{2})&=-\delta_{2}[B_{j,11}^{2}]\zeta +\delta_{2}^{2}([B_{j,11}^{2}][B_{j,12}^{2}]-[B_{j,11}B_{j,12}]^{2})\\
&<-\frac{1}{5}\zeta ^{\frac{35}{17}}|\widetilde{B}_{j}|^{-2}_{h'_{j}},
\end{array}
\end{equation*}
and hence
\begin{equation}\label{6.12}
  \sqrt{-d(\delta_{2})}>\frac{1}{\sqrt{5}}\zeta ^{\frac{35}{34}}|\widetilde{B}_{j}|^{-1}_{h'_{j}}.
\end{equation}
In view of Proposition 18 of \cite{P06},~there exists $\mathcal{P}\in SL(2,\mathbb{R})$,~with $\|\mathcal{P}\|\leq2(\frac{\|b_{0}-\delta_{2}b_{1}\|}{\sqrt{-d(\delta_{2})}})^{\frac{1}{2}}$
such that
$$\mathcal{P}^{-1}e^{b_{0}-\delta_{2}b_{1}}\mathcal{P}=\left(
                                                         \begin{array}{cc}
                                                           e^{ \sqrt{-d(\delta_{2})}} & 0 \\
                                                           0 & e^{- \sqrt{-d(\delta_{2})}} \\
                                                         \end{array}
                                                       \right).
$$
Since $|\widetilde{B}_{j}|^{6}_{h'_{j}}\zeta ^{\frac{1}{17}}\leq\frac{1}{40}$,~we have
$$\|b_{0}-\delta_{2}b_{1}\|\leq \zeta +\zeta ^{\frac{18}{17}}(1+\zeta )|\widetilde{B}_{j}|^{2}_{h'_{j}}\leq2\zeta ,$$
and then,~by \eqref{6.12},~one has
$$\frac{\|b_{0}-\delta_{2}b_{1}\|}{\sqrt{-d(\delta_{2})}}\leq\frac{\sqrt{5}\cdot2\zeta }{\zeta ^{\frac{35}{34}}|\widetilde{B}_{j}|^{-1}_{h'_{j}}}=2\sqrt{5}|\widetilde{B}_{j}|_{h'_{j}}\zeta ^{-\frac{1}{34}}.$$
By \eqref{6.3} of Lemma \ref{lem-6.1},~we have
\begin{equation*}
\begin{array}{ll}
\mathcal{P}^{-1}\delta_{2}^{2}\|P_{1}\|_{C^{0}}\mathcal{P}&\leq8\sqrt{5}|\widetilde{B}_{j}|_{h'_{j}}\zeta^{-\frac{1}{34}}\zeta^{\frac{36}{17}}
(61c_{\tau}^{2}\gamma^{-6}|\widetilde{B}_{j}|_{h'_{j}}^{4}+
\zeta^{-\frac{18}{17}}\zeta^{2}|\widetilde{B}_{j}|^{2}_{h'_{j}})\\
&\leq 680\sqrt{5}c_{\tau}^{2}\gamma^{-6}|\widetilde{B}_{j}|_{h'_{j}}^{5}\zeta^{\frac{34}{68}}\tilde{h}_j^{-2(4\tau+1)}\\
&\leq -d(\delta_{2}).
\end{array}
\end{equation*}
Under the condition \eqref{c},~we have $|I_{k}(\lambda v)|\geq \zeta ^{\frac{18}{17}}\geq |\lambda|^{2}e^{-\frac{13}{6}|k|^{2s}}$ with $k$ satisfies (\ref{k-2.7}-\ref{k-2.12}).
\end{pf}

\section{APPENDIX: PROOF OF PROPOSITION \ref{lem3.2}}

\begin{pf}
For $sl(2,\mathbb{R})$ and $su(1,1)$ are isomorphic via algebraic conjugation through some matrix.~ Using this property, we will prove the Proposition on the isomorphism group $su (1,1) $. We will discuss this in two ways.

\noindent (1)\textbf{Non-resonant case:}~when $0<|n|\leq N=\frac{2|\ln\epsilon|}{h-h_{+}},$~we have
\begin{equation}\label{3.7}
	|2\rho-\langle n,\alpha\rangle|\geq\epsilon^{\frac{1}{10}}.
\end{equation}
and $\alpha\in\ DC_{d}(\kappa,\tau),$~so
\begin{equation}\label{3.8}
	|\langle n,\alpha\rangle|\geq\frac{\kappa}{|n|^{\tau}}\geq\frac{\kappa}{|N|^{\tau}}\geq\epsilon^{\frac{1}{10}}.
\end{equation}
Define $$\Lambda_{n}:=\left\{f\in C^{\omega}_{\tau}(\mathbb{T}^{d},su(1,1))\Bigg|f(\theta)=\sum_{k\in\mathbb{Z},0<|k|\leq N}{\hat{f}(k)e^{i<k,\theta>}}\right\}.$$
Simple calculation shows that:~if $Y\in\Lambda_{n},$~then $$|A^{-1}Y(\theta+\alpha)A-Y(\theta)|_{h}\geq\epsilon^{\frac{1}{10}}|Y(\theta)|_{h}.$$
There is $Y\in\mathcal{B}_{h},~~F^{re}(\theta)\in\mathcal{B}^{re}_{h}(\epsilon^{\frac{3}{10}})$ so that $$e^{Y(\theta+\alpha)}(Ae^{F(\theta)})e^{-Y(\theta)}=Ae^{F^{re}(\theta)},$$ and $|Y|_{h}\leq\epsilon^{\frac{1}{2}},|F^{re}|_{h}\leq 2\epsilon.$
And we know the non-resonant case when $|k|\leq N$,~The only non-zero term of Fourier coefficient $\hat{F}(k)$ is that $\hat{F}(0).$ So for the pre-truncated portion there are:$$(\mathcal{T}_{N}F^{re})(\theta)=\hat{F}^{re}(0),\quad \Vert \hat{F}^{re}(0)\Vert\leq2\epsilon.$$
The truncated part is incorporated into the remainder by shrinking the analytic radius: $$|(\mathcal{R}_{N}F^{re})(\theta)|_{h^{'}}=|\sum_{|k|\leq N}\hat{F}^{re}(k)e^{i\langle k,~\theta\rangle}|_{h^{'}}\leq\frac{1}{2}\epsilon^{2}.$$
Again
\begin{equation}
\begin{array}{ll}
e^{F^{(re)}(\theta)} & =e^{\widehat{F}^{(re)}(0)+\mathcal{R}_{N}F^{(re)}(\theta)}\\
& =e^{\widehat{F}^{(re)}(0)}e^{F_{+}(\theta)}
\end{array}
\end{equation}	
so
$$|F_{+}(\theta)|_{h^{'}}\leq2|\mathcal{R}_{N}F^{(re)}(\theta)|_{h^{'}}\leq\epsilon^{2},$$
and $A_{+}=Ae^{\widehat{F}^{(re)}(0)}$ with the following estimates:
$$\|A-A_{+}\|\leq\|A\|\| I-e^{\widehat{F}^{(re)}(0)}\|\leq2\epsilon\|A\|.$$

\noindent (2)\textbf{Resonance case:} In this situation, only the case where $A$ is an ellipse needs to be considered, and the eigenvalues corresponding to this matrix are $e^{i\rho},e^{-i\rho},$~where $\rho\in\mathbb{R}\backslash\left\{0\right\}.$ If $\rho\in i\mathbb{R},$~the conclusion is established by (1).~Our selection of truncation $N$ ensures that when $0<|\breve{k}|<N$,~ there is and only one $\breve{k}$ that satisfies $|2\rho-\langle \breve{k},\alpha\rangle|<\epsilon^{\frac{1}{10}}$.

Firstly, diagonalize $A$.
By $|2\rho-\langle \breve{k},\alpha\rangle|<\epsilon^{\frac{1}{10}}, \epsilon\leq\frac{c}{\Vert A\Vert^{D}}(h-h^{'})^{D\tau},$~have $$|\ln\epsilon|^{\tau}\epsilon^{\frac{1}{10}}\leq\frac{\gamma(h-h^{'})}{2^{\tau+1}}.$$ Therefore, we have $$\frac{\gamma}{|\breve{k}|^{\tau}}\leq|\langle \breve{k},\alpha\rangle|\leq\epsilon^{\frac{1}{10}}+2|\rho|\leq\frac{\gamma}{2|\breve{k}|^{\tau}}+2|\rho|,$$
$$|\rho|\geq\frac{\gamma}{4|\breve{k}|^{\tau}}.$$

Therefore, there exists $U\in SU(1,1),\Vert U\Vert\leq\frac{2\Vert A\Vert}{|\rho|}\leq\frac{8\Vert A\Vert |\breve{k}|^{\tau}}{\gamma}
$ such that $$UAU^{-1}=\begin{pmatrix}
	e^{i\rho}  &   0  \\
	0          &   e^{-i\rho}\\
\end{pmatrix}=A^{'}$$
Let $G=UFU^{-1},$~by $$\Vert A\Vert|\breve{k}|^{\tau}\leq\Vert A\Vert|N|^{\tau}\leq\epsilon^{-\frac{1}{10}},$$ we have
\begin{eqnarray*}
&&\Vert U\Vert\leq\frac{8\Vert A\Vert N^{\tau}}{\gamma}\leq\frac{1}{2}\epsilon^{-\frac{1}{100}},\\
&&|G|_{h}\leq\Vert U\Vert^{2}|F|_{h}\leq\frac{1}{4}\epsilon^{\frac{49}{50}}:=\epsilon^{'}.
\end{eqnarray*}
Secondly, eliminate non resonant terms.
Define
\begin{eqnarray*}
&&\Theta_{1}:=\left\{k\in\mathbb{Z}^{d}:|\langle k,\alpha\rangle|\geq\epsilon^{\frac{1}{10}}\right\},\\
&&\Theta_{2}:=\left\{k\in\mathbb{Z}^{d}:2\rho-|\langle k,\alpha\rangle|\geq\epsilon^{\frac{1}{10}}\right\}.
\end{eqnarray*}
For $Cocycle ~(\alpha,A^{'}e^{G(\theta)}),$~by $\epsilon^{\frac{1}{10}}\geq13\Vert A^{'}\Vert^{2}(\epsilon^{'})^{\frac{1}{2}},$~we have
\begin{equation} \label{3.9}
e^{Y(\theta+\alpha)}(A^{'}e^{G(\theta)})e^{-Y(\theta)}=A^{'}e^{G^{re}(\theta)},
\end{equation}
and $|Y|_{h}\leq(\epsilon^{'})^{\frac{1}{2}},|G^{re}|_{h}\leq2\epsilon^{'}.$

For the remaining resonant structures, new truncations can be taken $\widetilde{N}=2^{-\frac{1}{\tau}}\gamma^{\frac{1}{\tau}}-N\ll N,$~at this point, the structure of $G^{re}$ is
\begin{equation}\label{3.10}
G^{re}(\theta)=\widehat{G}^{re}(0)+G^{re}_{A}(\theta)+G^{re}_{B}(\theta),
\end{equation}
where $$\widehat{G}^{re}(0)=\begin{pmatrix}
	i\widehat{t}(0)  & 0\\
	0&-i\widehat{t}(0)\\
\end{pmatrix},
G^{re}_{A}(\theta)=\begin{pmatrix}
0  &\widehat{ \vartheta}(\breve{k})e^{i\langle \breve{k},\theta\rangle}\\
	\widehat{ \vartheta}(\breve{k})e^{-i\langle \breve{k},\theta\rangle}&0\\
\end{pmatrix},$$
$$G^{re}_{B}(\theta)=\sum_{|k|>\widetilde{N}}\widehat{G}^{re}(k)e^{i\langle k,\theta\rangle}.$$

Finally, perform rotational conjugation on the resulting $Cocycle$.
Define $$Z(\theta)=\begin{pmatrix}
	e^{-\frac{\langle \breve{k},\theta\rangle}{2}i}&0\\
	0&e^{\frac{\langle \breve{k},\theta\rangle}{2}i}\\
\end{pmatrix},$$
obviously $Z(\theta)\in C^{\omega}_{h}(\mathbb{T}^{d},PSL(2,\mathbb{R}))$ meet with $$|Z(\theta)|_{h^{'}}\leq2e^{\frac{1}{2}|\breve{k}|h^{'}}\leq2e^{\frac{1}{2}Nh^{'}}\leq2\epsilon^{-\frac{h^{'}}{h-h^{'}}}.$$
Perform rotational conjugation on $Cocycle~(\alpha,A^{'}e^{G^{re}(\theta)})$ to obtain $$Z(\theta+\alpha)A^{'}e^{G^{re}(\theta)}Z^{-1}(\theta)=\widetilde{A}e^{\widetilde{G}(\theta)}.$$
due to $$Z(\theta+\alpha)A^{'}e^{G^{re}(\theta)}Z^{-1}(\theta):=Z(\theta+\alpha)A^{'}Z^{-1}(\theta)Z(\theta)e^{G^{re}(\theta)}Z^{-1}(\theta),$$
have $$\widetilde{A}=Z(\theta+\alpha)A^{'}Z^{-1}(\theta)=\begin{pmatrix}
	e^{i(\rho-\frac{\langle \breve{k},\theta\rangle}{2})}  &0\\
	0&	e^{-i(\rho-\frac{\langle\breve{k},\theta\rangle}{2})} \\
\end{pmatrix},$$
$$\widetilde{G}(\theta)=\begin{pmatrix}
	i\widehat{t}  & 0\\
	0&-i\widehat{t}\\
	\end{pmatrix}+\begin{pmatrix}
	0 &\vartheta(\breve{k})\\
	\overline{\widehat{\vartheta}(\breve{k})}&0\\
\end{pmatrix}+Z(\theta)e^{G^{re}_{B}(\theta)}Z^{-1}(\theta).$$
Convert $SU(1,1)$ to $SL(2,\mathbb{R}):$
\begin{equation*}
	\begin{aligned}
		\widetilde{A}^{'}   &=M^{-1}\widetilde{A}M,\\
		H               &=M^{-1}(G^{re}(0)+ZG^{re}_{A}(\theta)Z_{-1})M,\\
		F&=M^{-1}ZG^{re}_{B}(\theta)Z^{-1}M,\\
		B&=M^{-1}(Ze^{Y}U)M.\\
		\end{aligned}
\end{equation*}

In the above process, first convert $A$ into a diagonal type through $U$,~then use $e^{Y(\theta)}$ to eliminate the resonance term, and finally perform rotational conjugation.~Let $B=M^{-1}(Ze^{Y}U)M,$~have
\begin{eqnarray*}
\label{3.11}
&&B(\theta+\alpha)(Ae^{F(\theta)})B^{-1}(\theta)=\widetilde{A}^{'}e^{H+F^{'}},\\
\label{3.12}
&&|H|_{h}\leq |M^{-1}G^{re}(0)M|_{h}+|M^{-1}ZG^{re}_{A}(\theta)Z^{-1}M|_{h}\leq2\epsilon,\\
\label{3.13}
&&|F^{'}|_{h^{'}}\leq|M^{-1}ZG^{re}_{B}(\theta)Z^{-1}M|_{h^{'}}\leq\epsilon^{'},\\
\label{3.14}
&&|B|_{h^{'}}\leq|Z|_{h^{'}}|e^{Y(\theta)}|_{h^{'}}|U|_{h^{'}}\leq2\epsilon^{-\frac{h^{'}}{h-h^{'}}}\frac{1}{2}\epsilon^{-\frac{1}{100}}
=\epsilon^{-\frac{1}{100}-\frac{h^{'}}{h-h^{'}}}.
\end{eqnarray*}

Again
\begin{equation}\label{3.15}
e^{H+F^{'}(\theta)}=e^{H}+O(F^{'}(\theta))=e^{H}(I+e^{-H}\mathcal{O}(F^{'}(\theta)))=e^{H}e^{F_{+}(\theta)},
\end{equation}
record
\begin{equation}\label{3.16}
A_{+}=\widetilde{A}^{'}e^{H}=e^{A^{''}},e^{F_{+}(\theta)}=I+e^{-H}\mathcal{O}(F^{'}(\theta)),
\end{equation}
Set $H=M^{-1}\left(
              \begin{array}{cc}
                ia_{+} & b_{+} \\
                \overline{b}_{+} & -ia_{+} \\
              \end{array}\right)M$
then have
\begin{eqnarray*}
\label{3.17}
 && |a_{+}|^{2}+|b_{+}|^{2}\leq4|F|_{h},\quad |b_{+}-\widehat{g}(\breve{k})|\leq400\epsilon^{-\frac{3}{10}|F|_{h}^{2}e^{-|\breve{k}|h}},\\
\label{3.18}
&&|F_{+}|_{h^{'}}\leq2|F^{'}|_{h^{'}}\leq2\epsilon\ll\epsilon^{10},\\
\label{3.19}
&&\Vert A^{''}\Vert\leq2(|\rho-\frac{\langle \breve{k},\omega\rangle}{2}|+|g^{re}(0)|_{h}+|Zg^{re}_{A}(\theta)Z^{-1}|_{h})\leq 2\epsilon^{\frac{1}{10}}.
\end{eqnarray*}
\end{pf}

\section*{Acknowledgement}
X. Hou was partially  supported by NNSF of China (Grant
 11371019, 11671395) and Self-Determined Research Funds of Central China Normal University (CCNU19QN078).


\

\

\noindent
{\footnotesize Xuanji Hou\\
{\footnotesize School of Mathematics and Statistics, \\
and Key Laboratory of Mathematical Sciences, }\\
{\footnotesize Central China Normal University, Wuhan 430079,
China}\\
{\footnotesize hxj@mail.ccnu.edu.cn}\\

\noindent
{\footnotesize Li Zhang\\
{\footnotesize School of   Mathematics and Statistics,\\
and Key Laboratory of Mathematical Sciences, }\\
{\footnotesize Central China Normal University, Wuhan 430079,
China}\\
{\footnotesize lizhangmath@mails.ccnu.edu.cn}\\

\end{document}